\documentclass[12pt]{article}

\usepackage{amsmath}
\usepackage{amssymb}
\usepackage{amsthm}
\usepackage{latexsym}
\usepackage{color}
\usepackage{dsfont}
\usepackage{graphicx}
\usepackage{appendix}
\usepackage{lipsum}

\DeclareSymbolFont{calletters}{OMS}{cmsy}{m}{n}
\DeclareSymbolFontAlphabet{\mathcal}{calletters}

\def\be{\begin{eqnarray}}
\def\ee{\end{eqnarray}}

\def\b*{\begin{eqnarray*}}
\def\e*{\end{eqnarray*}}

\newtheorem{Theorem}{Theorem}[section]
\newtheorem{Definition}[Theorem]{Definition}
\newtheorem{Proposition}[Theorem]{Proposition}

\newtheorem{Assumption}[Theorem]{Assumption}

\newtheorem{Lemma}[Theorem]{Lemma}

\newtheorem{Remark}[Theorem]{Remark}
\newtheorem{Example}[Theorem]{Example}

\newtheoremstyle{named}{}{}{\itshape}{}{\bfseries}{.}{.5em}{\thmnote{#3's }#1}
\theoremstyle{named}

\makeatletter \@addtoreset{equation}{section}

\newcommand{\brak}[1]{\left(#1\right)}    
\newcommand{\crl}[1]{\left\{#1\right\}}   

\usepackage{color}

\def \D{\mathbb{D}}
\def \E{\mathbb{E}}
\def \F{\mathbb{F}}

\def \L{\mathbb{L}}

\def \N{\mathbb{N}}

\def \P{\mathbb{P}}
\def \Q{\mathbb{Q}}
\def \R{\mathbb{R}}

\def \T{\mathbb{T}}

\def \Mf{\mathbf{M}}

\def \Pf{\mathbf{P}}

\def \Tf{\mathbf{T}}

\def\Ac{{\cal A}}
\def\Bc{{\cal B}}
\def\Cc{{\cal C}}
\def\Dc{{\cal D}}

\def\Fc{{\cal F}}

\def\Hc{{\cal H}}

\def\Lc{{\cal L}}
\def\Mc{{\cal M}}

\def\Oc{{\cal O}}
\def\Pc{{\cal P}}

\def\Sc{{\cal S}}
\def\Tc{{\cal T}}
\def\Uc{{\cal U}}

\def\Wc{{\cal W}}
\def\Xc{{\cal X}}

\def \a {\alpha}
\def \bt {\beta}
\def \g {\gamma}
\def \d {\delta}
\def \eps {\varepsilon}
\def \z {\zeta}
\def \et {\eta}

\def \l {\lambda}
\def \m {\mu}
\def \n {\nu}
\def \r {\rho}
\def \s {\sigma}
\def \t {\tau}

\def \vp {\varphi}

\def \om {\omega}

\def \Gm {\Gamma}
\def \Lb {\Lambda}

\def \Dl {\Delta}
\def \Om {\Omega}

\def \sb {\subset}              
                    
\def \sbe {\subseteq}          
\def \spe {\supseteq}        

\def \pq {\preceq}

\def \lro {\longrightarrow}

\def \0 {\mathbf{0}}
\def \1{{\mathbf{1}}}

\def \Omo{\Om_+}

\def \Mco{\Mc_+}

\def \Omh{\hat{\Om}}
\def \htt{\hat{\tau}}
\def \hP{\hat{\P}}
\def \hp{\hat{\Pi}}
\def \hF{\hat{\F}}
\def \hh{\hat{H}}
\def \hH{\hat{\Hc}}
\def \hA{\hat{\Ac}}
\def \hD{\hat{\Dc}}

\def \Pft {\Pf^{\T}}
\def \Pfp {\Pf_{\pq}^{\T}}
\def \Pfo {\Pf_{\pq}^{\T,+}}
\def \bep {\boldsymbol{\varepsilon}}
\def \bl {\boldsymbol{\lambda}}
\def \bm {\boldsymbol{\mu}}
\def \bn {\boldsymbol{\nu}}
\def \pp {\mathsf{P}}
\def \dd {\mathsf{D}}

\def \ddd {\mathsf{D}_3}

\def \x {\times}

\def \wo {\stackrel{\Wc_1^{\T}}{\longrightarrow}}
\def \woi{\stackrel{\Wc_1}{\longrightarrow}}

\def \so {\stackrel{S}{\longrightarrow}}
\def \sto {\stackrel{S^{\ast}}{\longrightarrow}}

\def \dwo {\stackrel{\ast}{\Longrightarrow}_D}

\def \supp {\text{supp}}
\def \proof {{\textit{Proof.} }}

\def \prop {\textit{Proposition }}

\def \lems {\textit{Lemmas }}

\def \ass {\textit{Assumption }}
\def \asss {\textit{Assumptions }}

\newcommand{\rmi}{{\rm (i)$\>\>$}}

\newcommand{\rmii}{{\rm (ii)$\>\>$}}

\title{Tightness and duality of martingale transport \\
on the Skorokhod space
\thanks{We gratefully acknowledge the financial support of the ERC 321111 Rofirm, the ANR Isotace, and the Chairs Financial Risks (Risk Foundation, sponsored by Soci\'et\'e G\'en\'erale) and Finance and Sustainable Development (IEF sponsored by EDF and CA). We are grateful to Adam Jakubowski and Alexander Cox for fruitful discussions and comments.}}

\author{
	Gaoyue Guo\thanks{CMAP, Ecole Polytechnique, France. guo@cmap.polytechnique.fr}
	\and Xiaolu Tan\thanks{CEREMADE, University of Paris-Dauphine, France.
	tan@ceremade.dauphine.fr}
	\and Nizar Touzi\thanks{CMAP, Ecole Polytechnique, France. nizar.touzi@polytechnique.edu}
}
\date{\today}

\begin{document}
\bibliographystyle{plain}

\maketitle

\abstract{The \textit{martingale optimal transport} aims to optimally transfer a probability measure to another along the  class of martingales. This problem is mainly motivated by the robust superhedging of exotic derivatives in financial mathematics, which turns out to be the corresponding Kantorovich dual. In this paper we consider the continuous-time martingale transport on the Skorokhod space of c\`adl\`ag paths. Similar to the classical setting of optimal transport, we introduce different dual problems and establish the corresponding dualities by a crucial use of the $S-$topology and the dynamic programming principle\footnote{During the preparation of the final version of this paper, we knew from Mete Soner about the latest version of Dolinsky \& Soner \cite{DS2} which includes the finitely-many marginal constraints setting.}.

\vspace{3mm}

\noindent {\bf Key words.} $S-$\textit{topology}, \textit{dynamic programming principle}, \textit{robust superhedging}

\vspace{2mm}

\noindent {\bf AMS subject classification.} 65B05, 60B10, 60G44, 91G20


\section{Introduction}\label{sec:intro}

Initialed by the famous work of Monge and Kantorovich, the \textit{optimal transport problem} optimizes the cost of the  transfer of mass from one location to another. Namely, let  $\Pc(\R^d)$ be the space of probability measures on the Euclidean space $\R^d$. For any given measures $\m$, $\n\in\Pc(\R^d)$, let
\be\label{def:plan_op}
\Pc(\m,\n)&:=&\Big\{\P\in\Pc(\R^{d}\times \R^d): \P\circ X^{-1}=\m \mbox{ and } \P\circ Y^{-1}=\n\Big\},
\ee
where $(X,Y)$ denotes the canonical process on $\R^{d}\times \R^d$, i.e. $X(x,y)=x$ and $Y(x,y)=y$ for all $(x,y)\in\R^d\times\R^d$. Then the optimal transport problem consists in optimizing the expectation of some measurable function $\xi: \R^d\times \R^d\to\R$ among all probability measures in $\Pc(\m,\n)$. Various related issues are studied, e.g. the general duality theory and optimality results, we refer to Rachev \& R\"uschendorf \cite{RR} and Villani \cite{Villani} for a comprehensive account of the literature. 

Recently, a martingale optimal transport problem was introduced in Beiglb\"ock, Henry-Labord\`ere \& Penkner \cite{BHLP} in discrete-time (see Galichon, Henry-Labord\`ere \& Touzi \cite{GHLT} for the continuous-time case), where a maximization problem is considered over a subset $\Mc(\m,\n):=\big\{\P\in\Pc(\m,\n): \E^{\P}[Y|X]=X,~ \P\mbox{-a.s.}\big\}$ :
\b*
\pp(\m,\n)&:=&\sup_{\P\in\Mc(\m,\n)}\E^{\P}\big[\xi(X,Y)\big].
\e*
Each element of $\Mc(\m,\n)$ is called a \textit{transport plan}. Similarly to the classical setting, the corresponding dual problem is defined by
\b*
\dd(\m,\n)~:=~\inf_{(\l, \vp, H)\in\Dc(\m,\n)}\Big\{\int \l d\m+\int\vp d\n\Big\},
\e*
with $\Dc(\m,\n)$ being the collection of triplets $(\l, \vp, H)$, where $\l, \vp, H: \R^d\to\R$ are measurable functions such that $\l\in\L^1(\m)$, $\vp\in\L^1(\n)$ and 
\be\label{dualformulation}
\l(x)+\vp(y)+H(x)(y-x)~\ge~ \xi(x,y) \mbox{ for all } (x,y)\in\R^d\times\R^d.
\ee
The last dual formulation has the interpretation of minimal robust superhedging cost of derivative security defined by the payoff $\xi$ by trading the underlying security and any possible Vanilla option. When $d=1$, as observed by Breeden \& Litzenberger \cite{BL}, the marginal distributions of the underlying asset are recovered by the market prices of calls for all strikes, and any Vanilla option has a  non-ambiguous price as the integral of its payoff function with respect to the marginal.
Therefore, the inequality \eqref{dualformulation} represents a super-replication of $\xi$, which consist of the trading of the underlying and Vanilla options at different maturities. Since there is no specific model imposed on the process $(X,Y)$, the dual problem may be interpreted as the robust superhedging cost, i.e. the minimum cost to construct super-replications. Similar to the classical setting, the duality $\pp(\m,\n)=\dd(\m,\n)$ holds under quite general conditions.

The present paper considers the continuous-time martingale optimal transport problem. Let $\Xc:=\big\{\om=(\om_t)_{0\le t\le 1}: \om_t\in \R^d \mbox{ for all } t\in [0,1]\big\}$, where $\Xc$ is either the space of continuous functions or the Skorokhod space of c\`adl\`ag functions. Denote by $X=(X_{t})_{0\le t\le 1}$ the canonical process and by $\Mc$ the set of all martingale measures $\P$, i.e. $X$ is a martingale under $\P$. For a given family of probability measures $\bm=(\m_t)_{t\in\T}$, where $\T\sbe [0,1]$ is a subset, define by $\Mc(\bm)$ the subset of transport plans $\P$, i.e. $\P\circ X^{-1}_{t}=\m_t$ for all $t\in\T$. Then for a measurable function $\xi: \Xc\to\R$, the problem is defined by
\be\label{primalformulation}
\pp(\bm)&:=&\sup_{\P\in \Mc(\bm)}\E^{\P}[\xi(X)].
\ee
In contrast with the discrete-time case, the set $\Mc(\bm)$ is generally not tight with respect to the usual topologies. Without the crucial compactness, the arguments in the classical setting fail to be adapted to handle the related issues. 

In the existing literature, there are two dual formulations for the problem \eqref{primalformulation}, Galichon, Henry-Labord\`ere \& Touzi studied a class of transport plans defined by stochastic differential equations in \cite{GHLT} and introduced a quasi-sure dual problem. They applied a stochastic control approach and deduced the duality. Another important contribution is due to Dolinsky \& Soner \cite{DS, DS2}, see also Hou \& Obl\'oj \cite{Hou}, where the dual problem is still pathwisely formulated as in \eqref{dualformulation}. By discretizing the paths and a technical construction of approximated martingale measures, they avoid the compactness issue and derive the duality. 

In addition, the martingale optimal transport problem is studied by the approach of Skorokhod embedding problem. Following the seminal paper of Hobson \cite{Hobson}, this methodology generated developments in many directions, see e.g. Brown, Hobson \& Rogers \cite{BHR}, Cox \& Obl\'oj \cite{CO, CO1}, Cox, Hobson \& Obl\'oj \cite{CHO}, Cox, Obl\'oj \& Touzi \cite{COT}, Cox \& Wang \cite{CW}, Davis, Obl\'oj \& Raval \cite{DOR}, Gassiat, Oberhauser \& dos Reis \cite{GOR}, Hobson \& Klimmek \cite{HK1, HK2, HK3}, Hobson \& Neuberger \cite{HN} and Madan \& Yor \cite{MadanYor}. A thorough literature is provided in the survey papers Hobson \cite{Hobson1} and Obl\'oj \cite{Obloj}.

Our main contribution in the paper is to study systematically the tightness of the set $\Mc(\bm)$ by means of the $S-$topology introduced in Jakubowski \cite{Jakubowski}. Endowing properly the space of marginal laws with a Wasserstein kind topology, the tightness yields the upper semicontinuity of the map $\bm\mapsto \pp(\bm)$ and thus the first duality,  obtained by penalizing the marginal constraints. Based on the first duality and using respectively the dynamic programming principle and the discretization argument of path-space, the dualities are established for both quasi-sure and pathwise dual formulations.

The above analysis immediately gives rise to a stability consequence. Denote $\overline{\pp}:={\pp}$ and $\underline{\pp}(\bm):=\inf_{\P\in\Mc(\bm)}\E^{\P}[\xi(X)]$, then it is shown that the map $\bm\mapsto\overline{\pp}(\bm)$ (resp. $\bm\mapsto\underline{\pp}(\bm)$) is upper (resp. lower) semicontinuous, which yields the stability, i.e. for any sequence $(\bm^n)_{n\ge 1}$ convergent to $\bm$, there exists a sequence $(\eps_n)_{n\ge 1} \sbe \R_+$ convergent to zero such that 
\b*
\big[\underline{\pp}(\bm^n),~ \overline{\pp}(\bm^n)\big] &\sbe& \big[\underline{\pp}(\bm)-\eps_n,~ \overline{\pp}(\bm)+\eps_n\big] \mbox{ for all } n\ge 1,
\e*
i.e. the interval of model-free prices is stable with respect to the market. 

The paper is organized as follows. We formulate the martingale optimal transport problem and provide the dual problems in Section 2. In Section 3, the duality results are presented and we reduce the infinitely-many marginal constraints to the  finitely-many marginal constraints. In Sections 4, 5 we focus on the finitely-many marginal case and provide all related proofs.


\section{Martingale optimal transport}

For all $0\le s< t$, denote by $\D([s,t],\R^d)$ the space of c\`adl\`ag functions defined on $[s,t]$ taking values in $\R^d$. Let $\Om:=\D([0,1], \R^d)$ with generic element denoted by $\om$. Denote further by $X:=(X_t)_{0\le t\le 1}$ the canonical process, i.e. $X_t(\om)=\om_t$ and by $\F:=(\Fc_t)_{0\le t\le 1}$ its natural filtration, i.e. $\Fc_t=\s(X_u, u\le t)$. Let $\Pc:=\Pc(\Om,\Fc_1)$ be the set of probability measures on $\Om$. A probability measure $\P\in\Pc$ is called a martingale measure if the canonical process $X$ is a martingale under $\P$. Denote by $\Mc$ the collection of all martingale measures.

\subsection{Peacock and martingale optimal transport}

Let $\Pf:=\Pf(\R^d)$ be the space of all probability measures $\m$ on $\R^d$ with finite first moment. A pair $(\m,\n)\in\Pf\x\Pf$ is said to be increasing in convex ordering if  
\b*
\int_{\R^d}\l(x)\m(dx)~~:=~~\m(\l)&\le&\n(\l)~~:=~~\int_{\R^d}\l(x)\n(dx)
\e* 
holds for every convex function $\l:\R^d\to\R$. This relation is denoted by $\m\pq\n$. Let $\T\sbe [0,1]$ be some subset containing $1$ and define the $\T-$product of $\Pf$ by 
\b*
\Pft&:=&\Big\{\bm~:=~(\m_t)_{t\in\T}: \m_t\in \Pf \mbox{ for all } t\in\T\Big\}.
\e*

\begin{Definition}\label{def:pcoc}
A family of probability measures $\bm=(\m_t)_{t\in\T}\in\Pft$ is called a peacock ($\T-$peacock) if $\m_s\pq\m_t$ holds for all $s, t \in\T$ such that $s \le t$. A peacock $\bm$ is said to be c\`adl\`ag if the map $t\mapsto\m_t$ is c\`adl\`ag on $\T$ with respect to the weak convergence. Denote by $\Pfp$ the set of all c\`adl\`ag peacocks.
\end{Definition}

For each peacock $\bm\in\Pfp$, define the set of transport plans
	\be\label{def:motm}
		\Mc(\bm) &:=& \Big\{ \P \in \Mc: \P \circ X_t^{-1} = \m_t \mbox{ for all } t \in \T \Big\}.
	\ee
We may assume without loss of generality that $\T$ is closed under the \textit{lower limit topology}, i.e. the topology generated by all half-open intervals $[s,t)\sbe [0,1]$, see e.g. Steen \& Seebach \cite{SS}. Indeed, denote by $\bar{\T}$ the closure of $\T$ under the lower limit topology, then it follows that the law of $X_t$ for $t \in \bar \T$ is uniquely determined by the right continuity of $X$. This implies that $\Mc(\bar{\bm})=\Mc(\bm)$, where $\bar \bm:=(\bar \m_t)_{t\in\bar \T}$ is defined by
		\be\label{ex-pcoc}
			\bar \m_t~:=~\lim_{n\to\infty}\m_{t_n} \mbox{ for any sequence } (t_n)_{n\ge 1}\sbe\T \mbox{ decreasing to } t.
		\ee

\begin{Remark}\label{rem:Kellerer}
\rmi Since $\m_{t_n}\preceq \m_1$ for all $n$, we have
\b*
\m_{t_n}\big((x_i-K)^+\big) &\le& \m_{1}\big((x_i-K)^+\big) \mbox{ for all } i=1,\cdots, d,
\e*
thus showing that the sequence $(\m_{t_n})_{n\ge 1}$ is uniformly integrable. In particular, $(\m_{t_n})_{n\ge 1}$ is tight, and we may verify immediately by a direct density argument that any two possible accumulation points $\bar \m_t$ and $\bar \m_t'$ coincides, i.e. $\bar \m_t=\bar \m_t'$. Hence the sequence $(\m_{t_n})_{n\ge 1}$ converges weakly, justifying the convergence in \eqref{ex-pcoc} is well defined.

\vspace{1mm}

\noindent \rmii When $\T=[0,1]$, $\Mc(\bm)$ is nonempty by Kellerer's theorem, see e.g. Hirsch \& Roynette \cite{HR, HR2} and Kellerer \cite{Kellerer}. For a general closed $\T$, we may extend $\bm$ to some $\bar \bm = (\bar \m_t)_{0\le t\le1}$ by $\bar \m_t := \bar \mu_{\bar t}$ with $\bar t := \inf \{ s \ge t: s \in \T\}$. Clearly, $\bar \bm\in\Pf_{\pq}^{[0,1]}$ and $\bar \m_t=\m_t$ for all $t\in\T$. Hence $\Mc(\bm)\spe\Mc(\bar\bm)$ is again nonempty.

	\end{Remark}

Let $\xi: \Om\to\R$ be a measurable function. For every peacock $\bm\in \Pfp$, define the martingale optimal transport problem by
	\be\label{def:pp}
		\pp(\bm)&:=& \sup_{\P\in\Mc(\bm)}\E^{\P} \big[ \xi(X) \big],
	\ee
	where $\E^{\P}[\xi]:=\E^{\P}[\xi^+]-\E^{\P}[\xi^-]$ with the convention $+\infty-\infty=-\infty$.

\subsection{Dual problems}

\paragraph{First dual problem}

Let $\Lb$ be the set of continuous functions $\l:\R^d\to\R$ with linear growth, i.e. 
$\sup_{x\in\R^d} \big(|\l(x)|/(1+|x|) \big)<+\infty$. Define
\b*
\Lb^{\T}&:=&\Big\{\bl:=(\l_{t_i})_{1\le i\le m}:  t_i\in\T,~ \l_{t_i}\in\Lb \mbox{ for all } i=1,\cdots, m,~ m\in\N\Big\}.
\e*
For every $\bl=(\l_{t_i})_{1\le i\le m}\in\Lb^{\T}$, $\bm=(\m_t)_{t\in\T}\in\Pfp$ and $\om\in\Om$, denote
\b*
\bl(\om)~:=~\sum_{i=1}^m\l_{t_i}(\om_{t_i}) &\mbox{and}& \bm(\bl)~:=~\sum_{i=1}^m\m_{t_i}(\l_{t_i}).
\e*
Next, we introduce three dual formulations. Roughly speaking, as $X$ is required to be a martingale and has the given marginal laws in problem \eqref{def:pp}, then we dualize respectively these two constraints. The first dual problem is defined by
\be\label{def:dp0}
\dd_1(\bm)&:=&\inf_{\bl\in\Lb^{\T}}\Big\{\bm(\bl)+\sup_{\P\in\Mc} \E^{\P} \big[\xi(X) -\bl(X) \big] \Big\}.
\ee
The dual problem $\dd_1$ is the Kuhn-Tucker formulation in convex optimization, where the marginal constraints $\bm$ are penalized by the Lagrange multipliers $\bl$. 

\paragraph{Second dual problem}

The second dual problem dualizes further the martingale constraint and has close analogues in the mathematical finance literature in the context of a financial market with $d$ risky assets, where the price process is modeled by the canonical process $X=(X_t)_{0\le t\le 1}$. For technical reasons, the underlying process $X$ is assumed to be non-negative and start at some fixed price that may be normalized to be $\1:=(1, \cdots, 1)\in\R^d$. Namely, define the set of market scenarios 
\b*
\Omo&:=&\big\{\om\in\Om: \om_0=\1 \mbox{ and } \om_t\in\R^d_+ \mbox{ for all } t\in [0,1]\big\}
\e*
and the set of all possible models $\Mco:=\big\{\P\in\Mc: \supp(\P)\sbe\Omo\big\}$. Consequently, the market calibration $\bm$ should satisfy 
\be \label{con:muom+}
\m_0(dx)~=~\d_{\1}(dx)&\mbox{and}&\supp(\m_1)~\sbe~\R^d_+.
\ee
Moreover, let us denote by $\F^U = (\Fc^U_t)_{0 \le t \le 1}$ the universally completed filtration, i.e. $\Fc_t^U := \cap_{\P \in \Pc} \Fc_t^{\P}$, where $\Fc_t^{\P}$ is the completed $\sigma-$field of $\Fc_t$ under $\P$.

\begin{Definition}
A process $S=(S_t)_{0\le t\le 1}$ is called a $\Mco-$supermartingale if it is $\F^U-$adapted and is a $\P-$supermartingale for all $\P\in\Mco$. Denote by $\Sc$ the collection of all $\Mco-$supermartingales and by $\Sc_0\sbe\Sc$ the subset of processes starting at $0$. Denote further
\b*
\Dc_2(\xi)&:=&\Big\{(\bl,S)\in\Lb^{\T}\times\Sc_0: \bl(\om)+S_1(\om)\ge \xi(\om) \mbox{ for all } \om\in\Omo\Big\}.
\e*
\end{Definition} 

For a peacock $\bm\in\Pfp$ satisfying \eqref{con:muom+} , the second dual problem is defined by
\be\label{def:dp} 
\dd_2(\bm)&:=&\inf_{(\bl,S)\in\Dc_2(\xi)}\bm(\bl).
\ee
\begin{Remark}
\rmi Notice that the supermartingale $S\in\Sc$ is not required to have any regularity. If it were c\`adl\`ag then it would follow from Theorem 2.1 in Kramkov \cite{Kramkov} that, for every $\P\in\Mco$ there exist a predictable process $H^{\P}=(H^{\P}_t)_{0\le t\le 1}$ and an optional non-decreasing process $A^{\P}=(A^{\P}_t)_{0\le t\le 1}$ such that
\b*
S_t~=~S_0+\int_0^t H^{\P}_s dX_s-A^{\P}_t \mbox{ for all } t\in [0,1],~ \P\mbox{-a.s.}
\e*
However, it is not clear whether one can aggregate the last representation, i.e. find predictable processes $H$ and $A$ such that $(H,A)=(H^{\P},A^{\P})$, $\P$-almost surely. See also Nutz \cite{Nutz} for a partial result of this direction.

\vspace{1mm}

\noindent \rmii In financial mathematics, the pair $(\bl, H^{\P})$ has the interpretation of a semi-static super-replicating strategy under the model $\P$. If the aggregation above were possible, then the dual problem $\dd_2$ turns to the quasi-sure formulation of the robust superhedging problem, see also Beiglb\"ock, Nutz \& Touzi \cite{BNT}, and the duality $\pp=\dd_2$ reduces to the quasi-sure pricing-hedging duality.
\end{Remark}

\paragraph{Third dual problem}

Following the pioneering work \cite{Hobson} of Hobson, the martingale optimal transport approach is applied to study the robust hedging problems in finance. We do not postulate any specific model on the underlying assets and pursue here a robust approach. Assume further that all call/put options are liquid in the market for maturities $t\in\T$, thus yielding a family of marginal distributions $\bm=(\m_t)_{t\in\T}$ that is considered to be exogenous, see e.g. Breeden \& Litzenberger \cite{BL}. Then, the time $0$ market price of any derivative $\l(X_t)$ is given by $\m_t(\l)$. Hence, the cost of a static strategy $\bl\in\Lb^{\T}$ is $\bm(\bl)$. 

The return from a zero-initial cost dynamic trading, defined by a suitable process $H=(H_t)_{0\le t\le 1}$, is given by the stochastic integral $(H\cdot X)$ which we define similarly to Dolinsky \& Soner \cite{DS2}. We restrict $H: [0,1]\to \R^d$ to be left-continuous with bounded variation. Then, we may define the stochastic integral by integration by parts:
\be \label{eq:stoch_integral}
			(H\cdot X)_t
			~:=~
			H_t\cdot X_t - H_0\cdot X_0 - \int_0^t X_u\cdot dH_u \mbox{ for all } t\in [0,1],
		\ee 
where $\int_0^{t} X_u \cdot d H_u$ refers to the scalar Lebesgue-Stieltjes integration.

\begin{Definition}\label{def:super-replication}
An $\F-$adapted process $H: [0,1] \x \Omo\to\R^d$ is called a dynamic strategy if $t \mapsto H_t(\om)$ is left-continuous and of bounded variation for every $\om\in\Omo$ and $(H\cdot X)$ is a supermartingale under every $\P\in\Mco$. Let $\Ac$ be the set of all dynamic strategies and define the set of robust super-replications
		\b*
\Dc_3(\xi)&:=&\Big\{(\bl, H)\in\Lb^{\T}\x\Ac: \bl(\om)+(H\cdot X)_1(\om)\ge \xi(\om) \mbox{ for all } \om\in\Omo\Big\}.
		\e*
	\end{Definition}

For a peacock $\bm\in\Pfp$ satisfying \eqref{con:muom+}, the third dual problem is defined by
	\be\label{def:ph}
		\dd_3(\bm)
		&:=&
		\inf_{(\bl, H)\in\Dc_3(\xi)} \bm(\bl).
	\ee
	
\begin{Remark}
It is clear by definition that the weak duality $\pp(\bm)\le\dd_1(\bm)$ holds. Moreover, if the peacock $\bm$ satisfies \eqref{con:muom+}, then
\b*
\pp(\bm)~~\le~~\dd_1(\bm)~~\le~~\dd_2(\bm)~~\le~~\dd_3(\bm). 
\e*
\end{Remark}


\section{Main results}

	We aim to study the existence of optimal transport plans and establish the dualities in a systematic way. Before providing these results in Sections \ref{subsec:duality1} and \ref{subsec:duality2}, we first introduce some notions of topology on $\Om$ and the associated space of probability measures in Section \ref{subsec:Main_Prel}.

\subsection{Preliminaries}
\label{subsec:Main_Prel}

In the classical optimal transport problem, the relevant results (existence of optimizers, duality, etc.) rely essentially on the compactness condition of $\Mc(\m,\n)$. However, when passing to the  continuous-time case, as shown by Example \ref{example:SkoStopo} below, the set $\Mc(\bm)$ is in general not tight with respect to the topologies $L^{\infty}$ (uniform topology) and $J_1$ (Skorokhod topology). For our purpose, we endow $\Om$ with the $S-$topology introduced by Jakubowski \cite{Jakubowski} such that the Borel $\s$-field agrees with the projection $\s$-field $\Fc_1$, and more importantly, the $S-$topology facilitates the tightness issue and both Skorokhod representation theorem and Prohorov's theorem hold true. Before introducing the $S-$topology, we give an example which shows that the topologies $L^{\infty}$ and $J_1$ are not convenient to handle the tightness of $\Mc(\bm)$.

\begin{Example} \label{example:SkoStopo}
		Let $ M = (M_0, M_1,M_2)$ be a discrete-time martingale on some probability space such that 
		$\P \big[M_0 \neq M_1 \mbox{ and } M_1 \neq M_2\big] > 0$. 
		Define $\P_n:=\P\circ (M^n)^{-1}$ for $n\ge 3$, where $M^n=(M^n_t)_{0\le t\le 1}$ is defined by
		\b*
			M^n_t &:=& 
			M_0 \mathds{1}_{[0, \frac{1}{2} - \frac{1}{n})}(t)
			+
			M_1 \mathds{1}_{[ \frac{1}{2} - \frac{1}{n},  \frac{1}{2} + \frac{1}{n})} (t)
			+
			M_2 \mathds{1}_{[ \frac{1}{2} + \frac{1}{n}, 1]} (t) .
		\e*
Clearly, $\P_n\in\Mc(\bm)$ for all $n\ge 3$ with $\T=\{0,1\}$ and $\bm=\big(\P \circ M^{-1}_0,\P \circ M^{-1}_2\big)$. However, it follows from Theorem VI.3.21 in Jacod \& Shiryaev \cite{JacodShiryaev} that, the sequence $(\P_n)_{n\ge 3}$ is not $J_1-$tight and thus not $L^{\infty}-$tight.
\end{Example}

\begin{Definition}[$S-$topology] \label{def:Stopology}
		The $S-$topology on $\Om$ is the sequential topology induced by the following $S-$convergence, i.e. a set $F\sbe\Om$ is closed under $S-$topology if it contains all limits of its $S-$convergent subsequences, where the $S-$convergence (denoted by $\so$) is defined as follows.
		Let $(\om^n)_{n \ge 0}\sbe\Om$, we say that $\om^n\so\om^0$ as $n\to\infty$ if for each $\eps>0$, 
		we may find a sequence $(v^n_{\eps})_{n\ge 0} \sbe \Om$ such that
		\b*
			v^n_{\eps} \mbox{ has bounded variation, } \| \om^n - v^n_{\eps}\| ~\le~ \eps \mbox{ for all } n\ge 0
		\e*
		and
		\b*
			\lim_{n\to\infty}\int_{[0,1]} f(t) \cdot dv^n_{\eps}(t)~=~ \int_{[0,1]} f(t) \cdot dv^0_{\eps}(t) \mbox{ for all } f \in \Cc([0,1], \R^d).
		\e*
		We denote by $\sto$ the convergence induced by the $S-$topology.
\end{Definition}
	
\begin{Remark}
\rmi It is shown in Jakubowski \cite{Jakubowski} that the $S-$topology is not metrizable. However, its associated Borel $\sigma$-field coincides with $\Fc_1$. In a metric space, a subset is sequentially closed if and only if it is closed; 
	but in a non metrizable space, a sequentially closed set may not be closed.
	In particular, a sequentially closed set under $\sto$ may not be closed under $S-$topology (which is equivalent to be sequentially closed under $\so$).
	More precisely, it is shown in Remark 2.6 of \cite{Jakubowski} that the convergence $\sto$ is weaker than the original one $\so$. However, this is not a real problem for our case, since we know, from \cite{Jakubowski},
\b*
&&\om^n\sto\om, \mbox{ if and only if, in every subsequence } (n_k)_{k\ge 1}, \\
&&\mbox{one may find a further subsequence } (n_{k_l})_{l\ge 1} \mbox{ such that } \om^{n_{k_l}}\so\om.
\e*
In particular, a function $\xi:\Om\to\R$ is $S-$continuous (semicontinuous) if and only if $\xi$ is $S^{\ast}-$continuous (semicontinuous).

\vspace{1mm}

\noindent \rmii The functions $\om\mapsto\om_{i,1}$, $\om\mapsto\int_0^1\om_{i,t}dt$ and $\om\mapsto\int_0^1|\om_t|dt$ for $i=1,\cdots, d$ are $S-$continuous. The functions $\om\mapsto\|\om\|$ and $\om\mapsto\sup_{0\le t\le 1}\om_{i,t}$ for $i=1,\cdots, d$ are $S-$lower semicontinuous. 
\end{Remark}

Notice that the $S-$topology is not metrizable, then instead of the usual weak convergence, we use another convergence of probability measures introduced in \cite{Jakubowski}, which induces easy criteria for $S-$tightness and preserves the Prohorov's theorem, i.e. tightness yields sequential compactness.

\begin{Definition} \label{eq:CvgD}
	Let $(\P_n)_{n\ge 1}$ be a sequence of probability measures on the space $(\Om, \Fc_1)$. We say $\P_n\dwo \P$ if for each subsequence $(\P_{n_k})_{k\ge 1}$, one can find a further subsequence $(\P_{n_{k_l}})_{l\ge 1}$ and stochastic processes $( Y^l)_{l\ge 1}$ and $Y$ defined on the probability space $\big( [0,1], \Bc_{[0,1]}, \ell \big)$   endowed with the Lebesgue measure $\ell$, 
		such that $\Lc(Y^l)=\P_{n_{k_l}}$ for all $l\ge 1$, $\Lc(Y)=\P$, 
		\b*
			Y^l(e) ~\sto~Y(e) \mbox{ for all } e \in [0,1],
		\e*
		and for each $\eps>0$, there exists an $S^{\ast}-$compact subset $K_{\eps}\sbe\Om$ such that
		\b*
			\inf_{l\ge 1}\P_{n_{k_l}}\big[X\in K_{\eps}\big]&>& 1 - \eps.
		\e*
\end{Definition}

It follows from Jakubowski \cite{Jakubowski} (see Theorem \ref{Thm:Jakubowski}) that the convergence $\dwo$ implies in some sense the convergence of finite dimensional distributions that is specified later, and more importantly, the limit of every convergent sequence of martingale measures is still a martingale measure.

\begin{Remark}
	Meyer \& Zheng \cite{MeyerZheng} have also introduced a topology on $\Omega$ (for the case $d=1$), 
	called pseudo-path topology, 
	by considering the occupation measure induced by every path $\om \in \Om$ on $[0,1] \x \R$.
	We notice that $\om_n \to \om_0$ under $S-$topology induces $\om_n \to \om_0$ under the pseudo-path topology, and hence under the pseudo-path topology, it is easier to obtain the relative compactness of a sequence of martingale measures, but one has less continuous functionals defined on $\Om$.
	In particular, the simple maps $\om\mapsto\|\om\|$, $\om \mapsto \om_1$ are not upper semicontinuous,
	which makes it unsuitable to study the current martingale optimal transport.
\end{Remark}

We next introduce the Wasserstein distance for the purpose of deriving the duality $\pp=\dd_1$. Recall the set $\Pc(\m,\n)$ introduced in \eqref{def:plan_op}.

\begin{Definition}
The Wasserstein distance of order $1$ is defined by 
\b*
\Wc_1(\m,\n)~:=~\inf_{\P\in\Pc(\m,\n)}\E^{\P}\big[\big|X-Y\big|\big]
\mbox{ for all }
\m,~ \n\in\Pf.
\e*
A sequence $(\m^n)_{n\ge 1}\sbe\Pf$ converges to $\m\in\Pf$ in $\Wc_1$ if $\Wc_1(\m^n,\m)\to 0$ as $n\to\infty$ or, equivalently, $\lim_{n\to\infty}\m^n(\l)=\m(\l)$ for all $\l\in\Lb$, see e.g. Theorem 6.9 in Villani \cite{Villani}.
\end{Definition}

For $\big(\bm^n=(\m_t^n)_{t\in\T}\big)_{n \ge 1} \sbe \Pf_{\pq}^{\T}$ and $\bm=(\m_t)_{t\in\T}\in \Pf_{\pq}^{\T}$, we say that $\bm^n$ converges to $\bm$ if $\m^n_t$ converges to $\m_t$ in $\Wc_1$ for all $t\in\T$ and this convergence is denoted by $\wo$. We now provide a crucial tightness result for the present paper which is a consequence from \cite{Jakubowski}. 

Let $\T_0\sbe\T$ be the collection of all condensation points under the lower limit topology, i.e. $t=1$ or $[t,t+\eps)\cap\T$ is uncountable for any $\eps>0$.

	\begin{Lemma} \label{lemm:S_tight}
	Let $(\P_n)_{n \ge 1}$ be a sequence of probability measures such that $\P_n \in \Mc(\bm^n)$ for some $\bm^n\in\Pf_{\pq}^{\T}$, satisfying
	\be \label{eq:muW1}
			\bm^n~\wo~ \bm\in \Pf_{\preceq}^{\T}.
		\ee
	\rmi Then, $(\P_n)_{n \ge 1}$ is $S-$tight, i.e. any subsequence $(\P_{n_k})_{k \ge 1}$ admits a further  convergent subsequence under $\dwo$. Moreover, any limit point $\P$ of $(\P_n)_{n\ge 1}$ is again a martingale measure.
		
\vspace{1mm}

\noindent \rmii Assume in addition that $\T_0 = \T$, then $\P \in \Mc(\bm)$.
\end{Lemma}

\noindent \proof  \rmi By Theorem \ref{Thm:Jakubowski}, it is clear that $(\P_n)_{n \ge 1}$ is $S-$tight
	and there exist a convergent subsequence $(\P_{n_k})_{k \ge 1}$ with limit $\P\in\Pc$.
	Moreover, one has a countable subset $\Tc \sbe [0, 1)$
	such that for any finite set $\{u_1, \cdots, u_r \} \subset [0,1] \setminus \Tc$,
	\be \label{eq:cvg_finite_marginals}
		\P_{n_k} \circ \big(X_{u_1}, \cdots, X_{u_r} \big)^{-1}
		~\stackrel{\Lc}{\lro}~
		\P \circ \big(X_{u_1}, \cdots, X_{u_r} \big)^{-1} \mbox { as } k\to\infty.
	\ee
	Let $s, t \in [0,1] \setminus \Tc$ such that $s < t$,
	and take a finite subset $\{u_1, \cdots, u_r \}\sbe[0,s]\backslash\Tc$ 
	and a sequence of bounded continuous functions $\{ f_i\}_{1\le i\le r}$.
	Notice that for every $u \in [0,1]$,  $X_u$ is uniformly integrable with respect to $(\P_{n})_{n\ge 1}$. Indeeed, 
	\be \label{eq:unif_conv}
		\lim_{R \to \infty} 
		\sup_{n \ge 1} \E^{\P_n} \Big[ \big| X_u \big| 1_{|X_u| \ge R} \Big]
		&\le&
		\lim_{R \to \infty}
		\sup_{n \ge 1} \E^{\P_n} \Big[ \big( |X_u| - R/2 \big)_+ \Big] \nonumber\\
		&\le&
		\lim_{R \to \infty}
		\sup_{n \ge 1} \E^{\P_n} \Big[ \big( |X_1| - R/2 \big)_+ \Big] \nonumber\\
		&=&
		\lim_{R \to \infty}
		\sup_{n \ge 1} \m^n_1\big( (|x| - R/2)_+\big)		
		~~=~~ 0.
	\ee
 Combining \eqref{eq:cvg_finite_marginals} and \eqref{eq:unif_conv}, one has
	\b*
		\E^{\P} \Big[
			f_1(X_{u_1}) \cdots f_r(X_{u_r}) \big(X_t - X_s\big) 
		\Big]
		~=~
		\lim_{k\to\infty}
		\E^{\P_{n_k}} \Big[
			f_1(X_{u_1}) \cdots f_r(X_{u_r}) \big(X_t - X_s \big) 
		\Big]
		~=~
		0.
	\e*
	Since $\Tc$ is at most countable, it follows that $\E^{\P} [ X_t | \Fc_s] = X_s$ for any
	$s, t \in [0,1] \setminus \Tc$ such that $s < t$.
	It follows by the right continuity of $X$ that $\P \in \Mc$.

\vspace{1mm}

	\noindent \rmii To prove that $\P \in \Mc(\bm)$, it remains to show that $\P\circ X_t^{-1}=\m_t$ for all $t\in\T$.
	When $t \in \T \setminus \Tc$, 
	by the convergence \eqref{eq:cvg_finite_marginals} and the fact that $\bm^n\wo \bm$,
	it follows that $\P \circ X_t^{-1} = \mu_t$.
	Further, notice that $\T_0 = \T$, then for every $t \in \T$,
	there exists a sequence $(t_i)_{i\ge1}\sbe\T\backslash \Tc$ decreasing to $t$.
	Using again the right continuity of $X$, we conclude
	$\P \circ X_t^{-1} = \lim_{i \to \infty} \P \circ X^{-1}_{t_i} = \mu_t$.
	\qed
	
\vspace{2mm}	

\noindent As a consequence, the set $\Mc(\bm)$ is $S-$tight and it is closed if $\T_0 = \T$. The following example shows that the  closeness may fail when $\T_0 \neq \T$.
	
	\begin{Example} \label{example:closeness}
		Let $\T=\{0, 1\}$ and consider a random variable $Y$ such that $\P(Y= 1) = \P(Y = -1) =1/2$.
		Define $\P_n:=\P\circ (M^n)^{-1}$ for $n\ge 1$, where $M^n=(M^n_t)_{0\le t\le 1}$ is defined by
		\b*
			M^n_t &:=& Y\mathds{1}_{[\frac{1}{n}, 1]}(t).
		\e*
		Define a peacock $\bm=(\m_0,\m_1)$ by $\m_0:=\d_{\{0\}}$ 
		and $\m_1=\big(\d_{\{-1\}}+\d_{\{1\}}\big)/2$. Obviously, $\P_n\in\Mc(\bm)$ for all $n\ge 1$. 
		However, the limit of $(\P_n)_{n\ge 1}$ is a martingale measure $\P$ such that $X_t=X_0$, $\P$-a.s. and  $\P\circ X_0^{-1}=\mu_1$, which does not lie in $\Mc(\bm)$.
	\end{Example}


\subsection{Finitely-many marginal constraints}
\label{subsec:duality1}

We start by studying the finitely-many marginal case and assume throughout this subsection that $\T=\{0=t_0<\cdots<t_m=1\}$. Denote $\Dl t_i:=t_i-t_{i-1}$ for all $i=1,\cdots, m$ and $\Dl \T:=\min_{1\le i\le m}\Dl t_i$. Let us formulate some conditions on the reward function $\xi$. We shall see later that the usual examples satisfy our conditions.
	
\begin{Assumption}\label{ass:xi_S-usc}
$\limsup_{n\to\infty}\xi(\om^n)\le\xi(\om)$ holds for all $(\om^n)_{n\ge 1}\sbe\Om$ and $\om\in\Om$ such that 
	\b*
	\om^n~\sto~\om \mbox{ and } \om^n_{t_i}~\lro~\om_{t_i} \mbox{ for all } i=0,\cdots, m-1.
	\e*
\end{Assumption}

For $\bep=(\eps_1,\cdots, \eps_m)\in\R_+^m$ such that  $|\bep|<\Dl\T$, let $f_{\bep}$ (forward function) and $b_{\bep}$  (backward function) be two non-decreasing functions defined on $[0,1]$:
\be
			f_{\bep}(t)&:=&\sum_{i=1}^m
			\mathds{1}_{(t_{i-1}, t_{i}]}(t)\Big(t_{i-1}+\frac{\Dl t_i}{\Dl t_i-\eps_i}\big(t-t_{i-1}-\eps_i\big)^+\Big),   \label{forward} \\
			b_{\bep}(t)&:=&\sum_{i=1}^m
			\mathds{1}_{(t_{i-1}, t_{i}]}(t)\Big(t_{i}-\Big(\Dl t_i-\frac{\Dl t_i}{\Dl t_i-\eps_i}\big(t-t_{i-1}\big)\Big)^+\Big). \label{backward}
		\ee
			
		\begin{Assumption}\label{ass:xi_shift}
		There is a continuous function $\a : \R_+ \to \R_+$ with $\a(0)=0$ such that the following inequality holds for any $\bep\in\R^m_+$ satisfying  $|\bep|<\Dl\T$ 
\b*
\big|\xi(\om) -\xi \big(\om_{f_{\bep}} \big)\big|&\le &
			\a(|\bep|)\Big(1+\sum_{i=0}^m|\om_{t_i}|+\int_0^{1}|\om_t|dt\Big),~~~~~~~~~~~~~(1) \\
\big|\xi(\om) -\xi \big(\om_{b_{\bep}} \big)\big|&\le&
			\a(|\bep|)\Big(1+\sum_{i=0}^m|\om_{t_i}|+\int_0^{1}|\om_t|dt\Big),~~~~~~~~~~~~~(2)
\e*
where $\om_{f_{\bep}}$ (resp. $\om_{b_{\bep}}$) denotes the composition of $\om$ and $f_{\bep}$ (resp. $b_{\bep}$).
			\end{Assumption}

\begin{Theorem}\label{th:duality} Le $\xi$ be bounded from above and satisfies \asss \ref{ass:xi_S-usc} and \ref{ass:xi_shift} (1). Then for all $\bm \in \Pf_{\preceq}^{\T}$:

\noindent	\rmi  The duality $\pp(\bm)=\dd_1(\bm)$ holds.

\vspace{1mm}
		
		\noindent \rmii Assuming further that $\xi$ is bounded, the duality $\dd_1(\bm)=\dd_2(\bm)$ holds for all $\bm$ satisfying \eqref{con:muom+}.
	\end{Theorem}

To establish the duality $\dd_1(\bm)=\dd_3(\bm)$, we need more regularity conditions on $\xi$. Define a distance $\r_{\T}$ on $\Om$ by
	\be\label{def:sko}
		\r_{\T}(\om,\om')
		&:=&
		\sum_{i=1}^m\r_{[t_{i-1},t_i]}(\om,\om')+\Big|\int_0^{1}\big(\om_u-\om'_u\big)du\Big| \mbox{ for all } \om, \om'\in\Om,
	\ee
where $\r_{[s,t]}: \D([s,t],\R^d)\times \D([s,t],\R^d)\to \R_+$ dentoes the Skorokhod metric on the space $\D([s,t],\R^d)$. Clearly, $|\om_{t_i}-\om'_{t_i}|\le \r_{\T}(\om,\om')$ for all $\om,\om'\in\Om$ and $i=1,\cdots, m$. 

\begin{Assumption}\label{ass:xi_Skorpkhod}
$\xi$ is locally $\r_{\T}-$uniformly continuous, i.e. for every $R>0$, there exists a continuous increasing function $h_R: \R_+\to\R_+$ with $h(0)=0$, such that
\b*
|\xi(\om)-\xi(\om')|~\le~h_R\left(\rho_{\T}(\om,\om')\right) \mbox{ for all } \|\om\|,~ \|\om'\|\le R.
\e*
\end{Assumption}

	\begin{Theorem}\label{th:pricing-hedging duality}
		Let $\xi$ be bounded and $\bm \in \Pfp$ satisfying \eqref{con:muom+}. Then under \asss \ref{ass:xi_S-usc}, \ref{ass:xi_shift} and \ref{ass:xi_Skorpkhod}, the duality $\pp(\bm)=\dd_3(\bm)$ holds.

	\end{Theorem}

\begin{Remark}
Using the pathwise Doob's inequality in Acciaio, Beiglb\"ock, Penkner, Schachermayer \& Temme \cite{ABPST}, the boundeness condition in Theorem \ref{th:pricing-hedging duality} may be removed when $\m_1(|x|^p)<+\infty$ for some $p>1$, see also Dolinsky \& Soner \cite{DS2}.
\end{Remark}
	

\subsection{Infinitely-many marginal constraints}
\label{subsec:duality2}

Using approximation techniques, we then obtain some results for the martingale transport problem under infinitely-many marginal constraints.

\begin{Proposition}\label{th:approximation}
		Let $\xi$ be $S^{\ast}-$upper semicontinuous and bounded from above. For all $\bm\in\Pfp$:

\noindent \rmi Assume that there exists an increasing sequence of finite sets 
		$\{\T_n\}_{n\ge 1}$
		such that $1\in\T_n \subseteq \T$ for all $n \ge 1$ and $\cup_{n\ge 1}\T_n$ is dense in $\T$ under the lower limit topology. Then
		\b*
			\lim_{n\to\infty}\pp(\bm^n)
			&=&
			\pp(\bm) \mbox{ with } \bm^n:=(\m_t)_{t\in\T_n}.
		\e*
 \rmii Assume $\T_0=\T$, then there exists an optimal transport plan $\P^* \in \Mc(\bm)$, i.e.
		\be\label{eq:dp1}
			\pp(\bm)&=&\E^{\P^*}[\xi(X)].
		\ee
	\end{Proposition}

\noindent \proof \rmi It follows by the definition of $\bm^n$ that $\pp(\bm^n)$ is non-increasing with respect to $n$.
	Take a sequence $(\P_n)_{n \ge 1}$ such that $\P_n \in \Mc(\bm^n)$ 
	and
	\b*
		\pp(\bm)~~\le~~\lim_{n\to\infty}\pp(\bm^n)
		~~=~~
		\lim_{n \to \infty} \E^{\P_n}[\xi].
	\e*
	By Lemma \ref{lemm:S_tight} \rmi, there is a convergent subsequence 
	$(\P_{n_k})_{k \ge 1}$ with some limit $\P \in \Mc$. It follows by the same arguments in the proof of Lemma \ref{lemm:S_tight} that $\P\in\Mc(\bm)$ and, the upper semicontinuity of $\xi$ yields
	\b*
		\lim_{n\to\infty}\pp(\bm^n)
		~~=~~
		\lim_{k \to \infty} \E^{\P_{n_k}} [\xi]
		~~\le~~
		\E^{\P}[ \xi]~~\le~~\pp(\bm).
	\e*
 \rmii  Take a maximizing sequence $(\P_n)_{n \ge 1} \subseteq \Mc(\bm)$, then we may get a limit point $\P^{\ast}$ and by Lemma \ref{lemm:S_tight} \rmii, $\P^{\ast}$ is the required optimal transport plan.
	\qed	
	
	\vspace{2mm}
	
\noindent Consequently, we obtain immediately the dualities for general $\T$ through \prop \ref{th:approximation}.

	\begin{Theorem} \label{coro:duality}
		Let $\xi$ be $S^{\ast}-$upper semicontinuous and bounded from above and  $\bm\in\Pfp$, consider an increasing sequence of finite sets $\{\T_n\}_{n\ge 1}$ such that $\T_n \sbe \T$, $\cup_{n\ge 1}\T_n$ is dense in $\T$, and set $\bm^n:=(\mu_t)_{t\in\T_n}$.
		
	\noindent \rmi Assume that $\pp(\bm^n)=\dd_1(\bm^n)$ for all $n\ge 1$. Then $\pp(\bm)=\dd_1(\bm)$.
	
	\vspace{1mm}
	
	\noindent \rmii Assume further that $\bm$ satisfies \eqref{con:muom+} and $\dd_1(\bm^n)=\dd_2(\bm^n)=\dd_3(\bm^n)$ for all $n\ge 1$. Then $\dd_1(\bm)=\dd_2(\bm)=\dd_3(\bm)$.
	\end{Theorem}
	
\noindent \proof It is enough to show (i). Notice by definition that $\dd_1(\bm^n)\ge \dd_1(\bm)$ for all $n \ge 1$, then it follows by \prop \ref{th:approximation} \rmi that
	\b*
		\pp(\bm)
		~~=~~
		\lim_{n \to \infty} \pp(\bm^n)
		~~\ge~~
		\dd_1(\bm).
	\e*
	Then the proof is fulfilled by the weak duality $\pp(\bm) \le \dd_1(\bm)$. \qed
	
\begin{Remark}
In the present setting, the marginal constraint $\bm=(\mu_t)_{t\in\T}$ is given by a family of joint distributions $\mu_t$ on $\R^d$. If we replace the probability distribution $\mu_t$ by, either $d$ marginal distributions $(\mu^1_t, \cdots, \mu^d_t)$ on $\R$, or a joint distribution $\vec \mu_{\vec t}$ on $\R^{l \x d}$ for some $\vec t:=(t_1, \cdots, t_l)$ with $0\le t_1<\cdots<t_l\le 1$, then all the arguments still hold true and we can obtain similar duality results as in Theorems \ref{th:duality}, \ref{th:pricing-hedging duality} and \ref{coro:duality}.
\end{Remark}

\section{The dualities $\pp=\dd_1=\dd_2$}

In the following, we focus on the finite-marginal case, i.e. $\T=\{0=t_0<\cdots<t_m=1\}$ and start by proving the first duality. To prove the equality $\pp=\dd_1$, we shall apply the following well-known result from convex analysis. 		

\begin{Theorem}[Fenchel-Moreau]
Let $(E,\Sigma)$ be a Hausdorff locally convex space and $F: E \to \R$ be a concave and upper semicontinuous function. Then $F$ is equal to its biconjugate $F^{\ast\ast}$ which is defined by
\b*
F^{\ast\ast}(e)&:=&\inf_{e^{\ast}\in E^{\ast}}\Big\{\langle e, e^{\ast}\rangle+\sup_{e'\in E}\big(F(e')-\langle e', e^{\ast}\rangle\big)\Big\}
\e*
and $E^{\ast}$ denotes the dual space of $E$.
\end{Theorem}

Next we show that the map $\bm\mapsto \pp(\bm)$ is $\Wc_1^{\T}-$upper semicontinuous and concave and then identify its dual space to be $\Lb^{\T}$ by $\langle \bm, \bl\rangle=\bm(\bl)$.


\subsection{Space of signed measures on $\R^d$ and its dual space}

	Let $\Mf$ denote the space of all finite signed Borel measures $\mu$ on $\R^d$ satisfying
	\b*
		\int_{\R^d}\big(1+|x|\big)|\m|(dx)&<&+\infty.
	\e*
	It is clear that $\Mf$ is a linear vector space.
	We endow $\Mf$ with a topology (of Wasserstein kind) induced by the following convergence:
	Let $(\mu^n)_{n\ge 0}\sbe\Mf$ be a sequence of bounded signed measures,
	we say $\mu^n$ converges to $\m^0$ if 
	\b*
		\lim_{n\to\infty}\int_{\R^d}\l(x)\mu^n(dx) ~=~ \int_{\R^d}\l(x)\mu^0(dx)
		\mbox{ for all }
		\lambda \in \Lb.
	\e*
	Notice that the above topology restricted on the subspace $\Pf \sbe \Mf$ of probability measures 
	is exactly that induced by the Wasserstein distance. As for the space $\Mf_0$ of all finite signed Borel measures on $\R^d$ equipped with the weak convergence topology, it is well known that its dual space $\Mf_0^{\ast}$ can be identified as the space of all bounded continuous functions $\Lb_0$, see e.g. Deuschel \& Stroock \cite{DeuschelStroock}. The following lemma identifies the dual space of $\Mf$.

	\begin{Lemma}\label{lemm:M_dual_space}
		The space $\Mf$ is a Hausdorff locally convex space,
		and the duality relation
\b*
(\l, \m)~\in~\Lb\x\Mf&\longmapsto& \m(\l)
\e*
determines a representation of $\Mf^{\ast}$  as $\Lb$.
	\end{Lemma}

 The proof is almost the same as that of $\Mf_0^{\ast}=\Lb_0$. 
	For completeness, we provide a short proof in Appendix.
For the finite set $\T$, let us endow $\Mf^{\T}$ with the product topology and obviously, the dual space of $\Mf^{\T}$ is given by $\Lb^{\T}$.


\subsection{Proof of the duality $\pp=\dd_1$ }

In preparation for the first duality, we show first the upper semicontinuity of  $\bm \mapsto \pp(\bm)$ in the context of Theorem \ref{th:duality} \rmi. For $\bep=(\eps_1,\cdots, \eps_m)\in\R_+^m$ such that $|\bep|<\Dl \T$, we introduce
	\b*
		\Mc^{\bep}(\bm)
		&:=&
		\Big\{
			\P\in\Mc(\bm):
			X_t=X_{t_k} \mbox{ on } [t_k, t_k+\eps_k) \mbox{ for all } k=0,\cdots, m-1,~ \P\mbox{-a.s.}
		\Big\}
	\e*
	and
	\b*
		\pp^{\bep}(\bm)&:=&\sup_{\P\in\Mc^{\bep}(\bm)}\E^{\P} \big[\xi(X)\big]. 
	\e*
	
	\begin{Proposition}\label{prop:usc}
		Let $\xi$ be bounded from above and satisfying \asss \ref{ass:xi_S-usc} and \ref{ass:xi_shift} (1), then $\bm \mapsto \pp(\bm)$ is $\Wc_1^{\T}-$upper semicontinuous on $\Pfp$.
	\end{Proposition}
	
	\noindent \proof \rmi
	First notice $\pp^{\bep}(\bm)\le \pp(\bm)$ since $\Mc^{\bep}(\bm) \sbe \Mc(\bm)$. 
	Next, for each $\P\in\Mc(\bm)$, define $\P^{\bep}:=\P\circ X_{f_{\bep}}^{-1}$, where $f_{\bep}$ is defined in \eqref{forward}. 
	It is clear that $\P^{\bep}\in\Mc^{\bep}(\bm)$ and $\E^{\P^{\bep}}[\xi(X)]=\E^{\P}[\xi(X_{f_{\bep}})]$. 
	 It follows by \ass \ref{ass:xi_shift} (1),
	\b*
		\E^{\P} \big[\xi(X) \big]
		&\le& 
		\E^{\P} \big[\xi(X_{f_{\bep}}) \big] + \a(|\bep|)\Big(1+(m+2)\E^{\P} \big[|X_1| \big]\Big) \\
		&=&
		\E^{\P^{\bep}} \big[\xi(X) \big] + \a(|\bep|)\big( 1+(m+2) \mu_1(|x|)\big) \\
		&\le&
		\pp^{\bep}(\bm) + \a(|\bep|)\big( 1+(m+2) \mu_1(|x|)\big),
	\e*
	which implies that
	\b*
		\pp(\bm)
		&=&
		\inf_{0<|\bep|<\Dl\T}\Big(\pp^{\bep}(\bm)+  \a(|\bep|)\big( 1+(m+2) \mu_1(|x|)\big)\Big).
	\e*
	\rmii In order to prove that $\bm\mapsto\pp(\bm)$ is upper semicontinuous, it suffices to verify that $\bm\mapsto\pp^{\bep}(\bm)$ is upper semicontinuous. To see this, let $(\bm^n)_{n\ge 1}\sbe\Pfp$ be a sequence such that $\bm^n\wo\bm \in\Pfp$.
	By definition, we have a sequence $(\P_n)_{n \ge 1}$ such that 
	$\P_n \in \Mc^{\bep}(\bm^n)$ and
	\b*
		\limsup_{n\to\infty} \pp^{\bep}(\bm^n)
		&=&
		\limsup_{n\to\infty} \E^{\P_{n}} \big[\xi \big]. 
	\e*
	Then one may find a convergent subsequence $(\P_{n_k})_{k\ge 1}$ with limit $\P \in \Mc$. It follows by exactly  the same arguments as in Lemma \ref{lemm:S_tight} \rmii that $\P\in\Mc^{\bep}(\bm)$. Since $\xi$ is bounded from above, then it follows from Fatou's lemma that
	\b*
		\limsup_{n\to\infty} \E^{\P_{n}} [\xi]~~=~~\lim_{k\to\infty} ~\E^{\P_{n_k}} [\xi]
		~~\le~~
		\E^{\P}[ \xi]
		~~\le~~
		\pp^{\bep}(\bm),
	\e*
	which concludes the proof. \qed

\vspace{2mm}

Now we are ready to provide the first duality $\pp(\bm) = \dd_1(\bm)$. To apply Fenchel-Moreau theorem, we need to embed $\Pfp$ to a locally convex space. Recall that $\Mf$ is the space of all finite signed measures $\m$ such that
	\b*
		\int_{\R^d}\big(1+|x|\big)|\m|(dx)&<&+\infty,
	\e*
	and $\Mf^{\T}$ is its $\T-$product.
	We then extend the map $\pp$ from $\Pfp$ to $\Mf^{\T}$ by
	\b*
		\widetilde{\pp}(\bm)
		&:=&
		\begin{cases}
			 \pp(\bm), ~~&\mbox{ if } \bm\in\Pft_{\preceq}, \\
			-\infty, ~~& \mbox{ otherwise}.
		\end{cases}
	\e* 
 	
	\noindent \textit{Proof of Theorem \ref{th:duality} \rmi.} The concavity of the map $\bm \mapsto \pp(\bm)$ is immediate from its definition. Together with the upper semicontinuity of Proposition \ref{prop:usc}, we may directly verify that the extended map $\widetilde{\pp}$ is also $\Wc_1^{\T}-$upper semicontinuous and concave. Then, combining the Fenchel-Moreau theorem and Lemma \ref{lemm:M_dual_space}, it follows that for all $\bm\in \Mf^{\T}$,
	\b*
		\widetilde{\pp}(\bm)&=&\widetilde{\pp}^{\ast\ast}(\bm),
	\e*
	where $\widetilde{\pp}^{\ast\ast}$ denotes the biconjugae of $\widetilde{\pp}$.
	In particular, for $\bm\in\Pfp$ one has
	\b*
		&&
		\pp(\bm) ~~=~~ \widetilde{\pp}(\bm) ~~=~~ \widetilde{\pp}^{\ast\ast}(\bm)\\
		&=&
		\inf_{\bl\in\Lb^{\T}}\big\{\bm(\bl)-\widetilde{\pp}^{\ast}(\bl)\big\}
		~~=~~
		\inf_{\bl\in\Lb^{\T}}\Big\{\bm(\bl)-\inf_{\bn\in\Mf^{\T}}\big\{\bn(\bl)-\widetilde{\pp}(\bn)\big\}\Big\} \\
		&\ge&
		\inf_{\bl\in\Lb^{\T}}
		\Big\{\bm(\bl)
			+
			\sup_{\bn\in\Pfp} \Big\{\sup_{\P\in\Mc(\bn)} \E^{\P} \big[\xi - \bl(X) \big] \Big\}
		\Big\} \\
		&=&
		\inf_{\bl\in\Lb^{\T}}\Big\{\bm(\bl)+\sup_{\P\in\Mc} \E^{\P} \big[\xi -\bl(X) \big] \Big\}
		~~=~~
		\dd_1(\bm)
		~~\ge~~
		\pp(\bm),
	\e*
	which yields $\pp(\bm) = \dd_1(\bm)$.
	\qed


\subsection{Proof of the duality $\dd_1=\dd_2$}

For technical reasons, we need to restrict the static strategy $\bl$ to a smaller class of functions $\Lb_{lip}^{\T}$ defined by
	\b*
		\Lb_{lip}^{\T}
		&:=&
		\Big\{ \bl=(\l_{t_i})_{1\le i\le m} \in \Lb^{\T}: \mbox{each}~ \l_{t_i} ~\mbox{is boundedly supported and Lipschitz}
		\Big\}.
	\e*

	\begin{Proposition}\label{prop:dplip}
		Under the conditions of Theorem \ref{th:duality} \rmii one has
		\be
			\dd_1(\bm)
			&=&
			\inf_{\bl\in\Lb_{lip}^{\T}}\Big\{\bm(\bl)+\sup_{\P\in\Mco}\E^{\P} \big[\xi - \bl(X) \big]\Big\}. 
		\ee
	\end{Proposition}

	\noindent \proof Clearly, by the definition of $\dd_1$ and the fact that $\m_0=\d_{\1}(dx)$ and $\supp(\m_1)\sbe\R^d_+$, one obtains by interchanging $\inf$ and $\sup$ that
	\b*
		\dd_1(\bm) 
		&\ge&\inf_{\bl\in\Lb^{\T}}\sup_{\P\in\Mco}\Big\{\bm(\bl)+\E^{\P} \big[\xi - \bl(X) \big]\Big\} \\		
		&\ge&\sup_{\P\in\Mco}\inf_{\bl\in\Lb^{\T}}\Big\{\bm(\bl)+\E^{\P} \big[\xi - \bl(X) \big]\Big\} \\
		&=&\pp(\bm)~~=~~\dd_1(\bm),
	\e*
		by Theorem \ref{th:duality} (i). Hence
		\b*
		\dd_1(\bm) 
		&=&\inf_{\bl\in\Lb^{\T}}\Big\{\bm(\bl)+\sup_{\P\in\Mco}\E^{\P} \big[\xi - \bl(X) \big]\Big\}.
	\e*
	Next for every $\bl=(\l_{t_i})_{1\le i\le m}\in\Lb^{\T}$,
	there exists some constant $L>0$ such that for every $1\le i\le m$,
	\b*
		\l_{t_i}^L(x)
		~:=~
		\l_{t_i}(x) 
		-
		L(1+\mathbf{1}\cdot x) 
		\le 0 \mbox{ for all } x\in\R_+^d.
	\e*
	Denote $\bl^L := (\lambda^L_{t_i})_{1 \le i \le m}$, then for every martingale measure $\P \in \Mc_+$, 
	we have
	\b*
		\bm(\bl) + \E^{\P} \big[\xi - \bl(X) \big]
		&=&
		\bm(\bl^L) + \E^{\P} \big[\xi - \bl^L(X) \big].
	\e*
	Further, for each $R > 0$, let $\psi_{R}: \R^d \to [0,1]$ be some continuous function such that
	\b*
		\psi_{R}(x)=1 \mbox{ whenever } |x|\le R
		\mbox{ and }
		\psi_{R}(x)=0 \mbox{ whenever } |x|>R+1.
	\e*
	Let $\bl^{L,R}:=\big( \lambda^{L,R}_{t_i} \big)_{1\le i\le m}$ 
	with $\lambda^{L,R}_{t_i}(x) := \l^L_{t_i}(x)\psi_{R}(x) \ge \lambda^L_{t_i}(x)$,
	then
	\b*
		\sup_{\P\in\Mco}\E^{\P} \big [ \xi - \bl^{L,R}(X) \big]
		&\le&
		\sup_{\P\in\Mco}\E^{\P} \big [ \xi - \bl^{L}(X) \big].
	\e*
	On the other hand, for all $\P \in \Mc_+$ we have by monotone convergence theorem
	\b*
		\lim_{R \to \infty}
		\E^{\P} \big [ \xi - \bl^{L,R}(X) \big]
		&=&
		\E^{\P} \big [ \xi - \bl^L(X) \big].
	\e*
	Hence
	\b*
		\lim_{R \to \infty}
		\sup_{\P\in\Mco}\E^{\P} \big [ \xi - \bl^{L,R}(X) \big]
		&=&
		\sup_{\P\in\Mco}\E^{\P} \big [ \xi - \bl^L(X) \big].
	\e*
	It follows that
	\b*
		\lim_{R \to\infty} \Big(
			\bm(\bl^{L,R}) 
			+ 
			\sup_{\P\in\Mco} \E^{\P} \big[ \xi - \bl^{L,R}(X) \big]
		\Big)
		&=&
		\bm(\bl^L)+\sup_{\P\in\Mco}\E^{\P}[\xi - \bl^L(X)\big] \\
		&=&
		\bm(\bl )+\sup_{\P\in\Mco}\E^{\P}[\xi - \bl (X)\big].
	\e*	
	Finally, by a convolution argument each $\l_{t_i}^{L,R}$ can be approximated uniformly by some Lipschitz function that is also boundedly supported,
	which yields the required result.	 \qed
	
	\vspace{2mm}

\noindent For all $(\om,t)\in\Omo\times [0,1]$, denote by $\Bc^{sem}_{\om,t}\sbe\Pc$ the set of probability measures $\P$ such that $\P\big[X_s=\om_s \mbox{ for all } 0\le s\le t\big]=1$ and $(X_s)_{s\ge t}$ is a non-negative semimartingale under $\P$. Denote further 
\b*
\Mc^{loc}_{\om,t}&:=&\Big\{\P\in\Bc^{sem}_{\om,t}: (X_s)_{s\ge t} \mbox{ is a local martingale under } \P \Big\}.
\e*
Write in particular $\Bc^{sem}=\Bc^{sem}_{\om,0}$ and $\Mc^{loc}=\Mc^{loc}_{\om,0}$.
 Let $\z:\Om\to\R$ be a measurable function and put
\be\label{def:supermartingale}
V_t(\om)&:=&\sup_{\P\in\Mc^{loc}_{\om,t}}\E^{\P}\big[\z(X)\big].
\ee

Our objective now is to show that the process $(V_t)_{0\le t\le 1}$ is $\F^U-$adapted and that the dynamic programming principle holds. To achieve this, we use the related results in Neufeld \& Nutz \cite{NN, NN2}. Let $\P\in\Bc^{sem}$ be a semimartingale measure with the triplet $(B^{\P}, C^{\P}, \nu^{\P})$ of predictable semimartingale characteristics, see e.g. Chapter II of Jacod \& Shiryaev \cite{JacodShiryaev}. Notice that
\b*
\Mc^{loc}&=&\Big\{\P\in\Bc^{sem}: B^{\P}_t=0 \mbox{ for all } t\in [0,1]\Big\}.
\e*
By Theorem 2.5 in \cite{NN}, the map $\P\mapsto (B^{\P}, C^{\P}, \nu^{\P})$ is measurable, then it follows that $\Mc^{loc}$ is Borel. Moreover, by the same arguments we have the following lemma.

\begin{Lemma} \label{lem:sem1}
The set $\big\{(\om, t, \P)\in\Om\times [0,1]\times\Pc(\Om): \P\in\Mc^{loc}_{\om,t}\big\}$ is Borel. 
\end{Lemma}

By Theorem 2.1 in \cite{NN2}, we have the following lemma.

\begin{Lemma} \label{lem:sem2}
Let $\P\in\Mc^{loc}_{\om,t}$ and $\tau$ be an $\F-$stopping time taking values in $[t,1]$.

\noindent \rmi There is a family of conditional probability $(\P_{\om})_{\om \in \Om}$ of $\P$ with respect to $\Fc_{\tau}$ such that $\P_{\om} \in\Mc_{\om,\tau(\om)}^{loc}$ for $\P$-a.e. $\om \in \Om$.

\noindent \rmii Assume that there exists a family of probability measures $(\Q_{\om})_{\om\in\Omo}$ such that 
\b*
\Q_{\om}\in\Mc_{\om,\tau(\om)}^{loc} ~\mbox{for}~\P\mbox{-a.e.}~ \om \in \Om,~ \mbox{and the map } \om\mapsto\Q_{\om} \mbox{ is}~ \Fc_{\tau}-\mbox{measurable},   
\e*
then $\P\otimes\Q\in\Mc^{loc}_{\om, t}$, where 
\b*
\P\otimes\Q(\cdot)&:=&\int_{\Om}\Q_{\om}(\cdot) \P(d \om).
\e*
\end{Lemma}

The dynamic programming principle follows by Lemmas \ref{lem:sem1} and \ref{lem:sem2}, and as a consequence we have the following proposition. 

\begin{Proposition}\label{lem:dyn-pro}
Assume that $\zeta$ is bounded, then the process $V=(V_t)_{0\le t\le 1}$ defined in \eqref{def:supermartingale} is a $\Mco-$supermartingale, i.e. $V\in\Sc$.
\end{Proposition}

\begin{Proposition}  \label{prop:dyn-pro}
Let $\zeta$ be a measurable and bounded function, then one has
\b*
\sup_{\P\in\Mco}\E^{\P}[\zeta(X)]&=&\inf\Big\{V_0: (V_t)_{0\le t\le 1}\in\Sc \mbox{ such that } V_1(\om)\ge \zeta(\om) \mbox{ for all } \om\in\Omo\Big\}.
\e*
\end{Proposition}

\noindent \proof By Proposition \ref{lem:dyn-pro} with the process $V$ defined in \eqref{def:supermartingale}, it remains to show that 
\b*
\sup_{\P\in\Mco}\E^{\P}[\zeta(X)]&=&\sup_{\P\in\Mc^{loc}}\E^{\P}[\zeta(X)].
\e*
It is clear that $\sup_{\P\in\Mco}\E^{\P}[\zeta(X)]\le\sup_{\P\in\Mc^{loc}}\E^{\P}[\zeta(X)]$ since $\Mco\sbe\Mc^{loc}$, then it suffices to prove the converse inequality. For each $\P\in\Mc^{loc}$, there exists an increasing sequence of stopping times $(\sigma_n)_{n\ge 1}$ such that $\sigma_n\to+\infty$, $\P$-almost surely and $X_{\sigma_n \wedge \cdot}$ is a $\P-$martingale, where $X_{\sigma_n\wedge\cdot}:=(X_{\sigma_n\wedge t})_{0\le t\le 1}$. Hence
\b*
\E^{\P}[\zeta(X_{\sigma_n\wedge\cdot})]&\le&\sup_{\Q\in\Mco}\E^{\Q}[\zeta(X)].
\e*
The required result follows from the dominated convergence theorem. \qed
	\vspace{3mm}

\noindent \textit{Proof of Theorem \ref{th:duality} \rmii.} It remains to show $\dd_1(\bm)\ge \dd_2(\bm)$. Indeed, one has by Proposition \ref{prop:dplip},
\b*
\dd_1(\bm)&=&\inf_{\bl\in\Lb_{lip}^{\T}}\Big\{\bm(\bl)+\sup_{\P\in\Mco} \E^{\P} \big[\xi(X) -\bl(X) \big] \Big\}.
\e*
For each $\eps>0$, by Proposition \ref{prop:dyn-pro} there exist a vector $\bl^{\eps}\in\Lambda_{lip}^{\T}$ and a process $V^{\eps}=(V^{\eps}_t)_{0\le t\le 1}\in\Sc$ such that
\b*
\dd_1(\bm)+\eps~\ge~ \bm(\bl^{\eps})+V^{\eps}_0 &\mbox{and}& V_1^{\eps}(\om)~\ge~ \xi(\om)-\bl^{\eps}(\om).
\e*
This implies that $\dd_1(\bm)+\eps\ge\dd_2(\bm)$, and the required result by the arbitrariness of $\eps$. \qed


\section{Proof of the duality $\dd_1=\dd_3$}

Now let us turn to prove the third duality $\dd_1= \dd_3$ in Theorem \ref{th:pricing-hedging duality}. We will follow the idea in Dolinsky \& Soner \cite{DS2} to discretize the underlying paths and then use the classical
constrained duality result of F\"ollmer \& Kramkov \cite{FK}. The proof in \cite{DS2} relies on the min-max theorem and the explicit approximation of a martingale measure. We emphasize that the present proof is less technically involved than \cite{DS2} as the marginals constraints have already been reduced by the first duality.


\subsection{Reduction of $\xi$ to be boundedly supported}

In this section we denote $\pp(\bm,\xi)$ and $\ddd(\bm,\xi)$ in place of $\pp(\bm)$ and $\ddd(\bm)$
to emphasize the dependence on $\xi$, then clearly for any $\xi$, $\xi': \Om\to \R$ and $c\in\R$, one has
\b*
\ddd(\bm, \xi+\xi')~\le~\ddd(\bm, \xi)+\ddd(\bm, \xi') &\mbox{and}& \ddd(\bm, \xi+c)~ =~\ddd(\bm, \xi)+c.
\e*
In particular for $c>0$ one has 
\b*
\ddd(\bm, c\xi)~=~c\ddd(\bm, \xi).
\e*
Hence, under the conditions of Theorem \ref{th:pricing-hedging duality}, we may assume that $0\le\xi\le 1$. Indeed, we show next that it suffices to establish the duality $\pp(\bm, \xi) =\ddd(\bm, \xi)$ for $\xi$ that is boundedly supported. For all $R>0$, define 
the continuous function $\chi_R: \R_+ \to [0,1]$ by
	\b*
		\chi_R(x)
		~:=~
		\mathds 1_{[0, R]}(x) + (R+1-x) \mathds 1_{(R, R+1]}(x) \mbox{ for all } x \in \R_+.
	\e*
	Denote further for $R>0$
	\b*
		\xi_R(\om)
		~:=~
		\xi(\om) \chi_R(\|\om\|) \mbox{ for all } \om \in \Om.
	\e*
Notice that $0\le \xi\le 1$ yields $\xi_R(\om) \le \xi(\om) \le \xi_R(\om) + \mathds 1_{\{ \|\om\| \ge R\}}$, then it follows that
	\be \label{eq:xiR_immed}
		\dd_3(\bm,\xi_R)
		~\le~
		\dd_3(\bm,\xi)
		~\le~
		\dd_3(\bm,\xi_R)
		+
		\dd_3(\bm,\mathds{1}_{\{\|X\|\ge R\}}). 
	\ee

	\begin{Lemma} \label{lemm:xi_to_xiR}
		Let $\xi$ be bounded and $\bm \in \Pfo$. Then
		\b*
			\dd_3(\bm,\xi)&=&\lim_{R\to\infty}\dd_3(\bm,\xi_R)
		\e*
	\end{Lemma}

	\noindent \proof  It is enough to prove by \eqref{eq:xiR_immed} that
	\b*
		\lim_{R\to+\infty} \dd_3(\bm,\mathds{1}_{\{ \|X\|\ge R\}}) 
		&=& 0.
	\e*
	This is indeed a direct consequence of the pathwise inequality, see e.g. Lemma 2.3 of Brown, Hobson and Rogers \cite{BHR}
		\b*
		\mathds 1_{\{ \|X_{i}\| \ge R\}}
		~\le~
		\frac{(|X_{i,1}| -K)^+}{R-K}
		+
		\mathds 1_{\{ \|X_{i}\| \ge R\}} \frac{R- X_{i,1} }{R-K} \mbox{ for all } i=1,\cdots, d
	\e*
holds for every $0< K < R$.
	It follows by taking $K=R/2d$ that
	\b*
		\dd_3(\bm,\mathds{1}_{\{\|X\| \ge R\}})~~\le~~\sum_{i=1}^d\dd_3(\bm,\mathds{1}_{\{\|X_i\| \ge R/d\}})
		~~\le~~
		\frac{2d}{R}\sum_{i=1}^d 
		\m_1\big((x_i-\frac{R}{2d})^+\big).
	\e*
	The proof is fulfilled by letting $R \to +\infty$.
	\qed

\vspace{2mm}

\noindent Next we show that $\xi_R$ inherits almost the same properties as $\xi$.
	
\begin{Lemma}\label{lem:zeta_ass}
For each $R>0$:

\noindent\rmi If $\xi$ satisfies \asss \ref{ass:xi_S-usc} and \ref{ass:xi_shift}, then so does $\xi_R$.
		
		\noindent \rmii If $\xi$ satifies \ass \ref{ass:xi_Skorpkhod}, then 
		\begin{equation}\label{ass:xi_sko}
		\begin{array}{c}
		\xi_R~~\mbox{is $L^{\infty}-$uniformly continuous, and}
		\\ 
		\xi_R(\om)-\xi_R(\om')~\le~\bt\big(\r_{\T}(\om,\om')\big) \mbox{ for all } \om,\om'\in\Om \mbox{ such that } \|\om'\|\le \|\om\|
		\end{array}
		\end{equation}	
for some continuous increasing function $\bt:\R_+\to\R_+$ with $\bt(0)=0$.
\end{Lemma}
	
	\noindent \proof
	\rmi follows by the fact that $\om\mapsto\|\om\|$ is $S^{\ast}-$lower semicontinuous and $\|\om_{f_{\eps}}\|=\|\om_{b_{\eps}}\|=\|\om\|$. Let us turn to show \rmii. Notice that $\xi$ is $\r_{\T}-$uniformly continuous on $\big\{\om: \|\om\|\le R\big\}$, i.e. there exists a continuous increasing function $\bt:\R_+\to\R_+$ with $\bt(0)=0$ such that for all $\|\om\|,\|\om'\|\le R$
\b*
|\xi(\om)-\xi(\om')|&\le&\bt\big(\r_{\T}(\om,\om')\big).
\e*
Hence, for any $\om,\om'\in\Om$ such that $\|\om'\|\le \|\om\|$, one has 
	\b*
	\xi_R(\om)-\xi_R(\om')~~\le~~\mathds{1}_{\{\|\om\|\le R\}}\big(\xi(\om)-\xi( \om')\big)~~\le~~\bt\big(\r_{\T}(\om,\om')\big).
	\e* 
Moreover,
\b*
|\xi_R(\om')-\xi_R(\om)|&\le & |\xi(\om')-\xi(\om)|\chi_R(\|\om\|)+|\chi_R(\|\om\|)-\chi_R(\|\om'\|)| \\
&\le &\bt\big(\r_{\T}(\om,\om')\big)+\|\om-\om'\| \\
&\le &\bt\big(2\|\om-\om'\|\big)+\|\om-\om'\|,
\e*	
which yields the $L^{\infty}-$uniform continuity of $\xi_R$.\qed

\vspace{2mm}

\noindent Therefore, in the following it suffices to consider the function $\xi$ that is boundedly supported such that the  Assumptions \ref{ass:xi_S-usc}, \ref{ass:xi_shift} and Condition \eqref{ass:xi_sko} hold. Similar to the proof of the duality $\pp (\bm) =\dd_2(\bm)$, it remains to prove a duality without marginal constraints.


\subsection{Duality without marginal constraints}

We consider in this section the optimization problem without marginal constraints. Let $\zeta:\Om\to\R$ be  bounded and define
	\be \label{eq:PD_delta}
		\pp(\zeta)
		~:=~
		\sup_{\P\in\Mco}\E^{\P} [ \zeta(X)]
		&\mbox{and}&
		\dd(\zeta)
		~:=~
		\inf_{(z, H) \in \Dc(\zeta)} z.
	\ee
	where, with the same definition of integral in \eqref{eq:stoch_integral},
	\b*
		\Dc(\zeta)
		&:=&
		\Big\{
			(z,H)\in \R\x\Ac
			: z+(H\cdot X)_1(\om)\ge \zeta(\om) \mbox{ for all } \om\in\Omo
		\Big\}.
	\e*
	We provide immediately a duality result for the above optimization problems,
	and leave its proof in Section \ref{subset:proof_duality_no_marginal}.

	\begin{Theorem}\label{th:dwc}
		Suppose that $\zeta$ satisfies the Assumptions \ref{ass:xi_S-usc}, \ref{ass:xi_shift} and Condition \eqref{ass:xi_sko}, then
		\be\label{eq:dwc}
			\pp (\zeta) &=& \dd(\zeta).
		\ee
	\end{Theorem}
			
By exactly the same arguments as in the proof of Theorem \ref{th:duality} \rmii, the duality $\pp(\bm,\xi_R)=\dd_3(\bm,\xi_R)$ follows immediately by taking $\zeta=\xi_R-\bl$ in Theorem \ref{th:dwc}.

	\vspace{3mm}

\noindent \textit{Proof of Theorem \ref{th:pricing-hedging duality}.} 
	Using Lemma \ref{lemm:xi_to_xiR} as well as the first duality $\pp=\dd_1$ for $\xi_R$, one has 
	\b*
	\pp(\bm,\xi)~~\ge~~\lim_{R\to\infty}\pp(\bm,\xi_R)
	~~=~~\lim_{R\to\infty}\dd_3(\bm,\xi_R)~~=~~\dd_3(\bm,\xi).
	\e*
Hence we conclude the proof by the weak duality $\pp(\bm,\xi)\le \dd_3(\bm,\xi)$. \qed


\subsection{Proof of Theorem \ref{th:dwc}}
\label{subset:proof_duality_no_marginal}

Recall that $\T=\{0=t_0<\cdots<t_m=1\}$, $\Dl t_i=t_i-t_{i-1}$ for $i=1,\cdots, m$ and $\Dl \T=\min_{1\le i\le m}\Dl t_i$. Let $\z:\Om\to\R$ be measurable and boundedly supported. Then for each $0\le \d<\Dl\T$, denote $\Om^{\d}:=\D([0,1+\d],\R^d)$ and all its elements by $\om^{\d}$. Put $\T^{\d}:=\{0=t_0^{\d}<\cdots<t_m^{\d}=1+\d\}$, where
$t_i^{\d}~:=~k_{\d}t_i$ for all $i=0,\cdots, m$ with $k_{\d}~:=~1+\d$. Define $\z^{\d}:\Om^{\d}\to\R$ by
\be\label{def:extension}
\z^{\d}(\om^{\d})~:=~\z(\bar\om^{\d}), \mbox{ where } \bar\om^{\d}\in\Om \mbox{ is defined by } \bar\om^{\d}_t~:=~\om^{\d}_{k_{\d}t} \mbox{ for all } t\in [0,1]. 
\ee

\begin{Proposition}\label{prop:extension}
Assume that $\z$ satisfies the \asss \ref{ass:xi_S-usc}, \ref{ass:xi_shift} and Condition \eqref{ass:xi_sko}. Then:

\noindent \rmi For all $0\le \d<\Dl\T$, the $\z^{\d}$ defined by \eqref{def:extension} satisfies the \asss \ref{ass:xi_S-usc}, \ref{ass:xi_shift} and Condition \eqref{ass:xi_sko}.

\noindent \rmii There is a continuous function $\et: \R_+\to\R_+$ with $\et(0)=0$ such that for all $0\le \d< \d'<\Dl \T$ the following inequality holds
\b*
\Big|\z^{\d}(\om^{\d})-\z^{\d'}(\om^{\d',\d})\Big|~\le ~ \et\Big(\frac{\d'-\d}{1+\d'}\Big)\Big(1+\sum_{i=0}^m|\om^{\d}_{t_i^{\d}}|+\int_0^{1+\d}|\om^{\d}_t|dt\Big) \mbox{ for all } \om^{\d}\in\Om^{\d},
\e*
where $\om^{\d',\d}\in\Om^{\d'}$ is defined by
\b*
\om^{\d',\d}_t~:=~\om^{\d}_{(t-t_i^{\d'}+t_i^{\d})\wedge t_{i+1}^{\d}} \mbox{ for all } t\in [t_i^{\d'}, t_{i+1}^{\d'}] \mbox{ and } i=0,\cdots, m-1.
\e*
\end{Proposition}

\noindent \proof (i) will be proved in Lemmas \ref{lem:ex1},  \ref{lem:ex2} and  \ref{lem:ex3} in Section \ref{subsec:extension}. 

\vspace{1mm}

\noindent \rmii Clearly, $\z^{\d'}(\om^{\d',\d})~=~\z(\bar \om^{\d',\d})$,
where
\b*
\bar \om^{\d',\d}_t~:=~\om^{\d',\d}_{k_{\d'}t} \mbox{ for all } t\in [0,1].
\e*
Direct computation reveals that 
 $\bar \om^{\d,'\d}=\bar \om^{\d}\circ b_{\bep}$ with
\b*
\bep&:=&\frac{\d'-\d}{1+\d'}\big(\Dl t_1,\cdots, \Dl t_m\big).
\e*
Hence by \ass \ref{ass:xi_shift} one obtains
\b*
\Big|\z^{\d}(\om^{\d})-\z^{\d'}(\om^{\d',\d})\Big|&=&\Big|\z(\bar \om^{\d})-\z(\bar \om^{\d',\d})\Big| \\
&\le&	\a(|\bep|)\Big(1+\sum_{i=0}^m|\bar\om^{\d}_{t_i}|+\int_0^{1}|\bar \om^{\d}_t|dt\Big) \\
  &\le&	\a\Big(|\Dl\T|\frac{\d'-\d}{1+\d'}\Big)\Big(1+\sum_{i=0}^m|\om^{\d}_{t_i^{\d}}|+\int_0^{1+{\d}}|\om^{\d}_t|dt\Big).
\e*
The proof is completed by taking $\et(\cdot)=\a(\Dl\T \times \cdot)$.  \qed

\vspace{2mm}

\noindent We are now ready to prove the required duality. Define
	 \b*
		\Omo^{\d}
		&:=&
		\big\{\om^{\d}\in\Om^{\d}: \om^{\d}_0=\1 \mbox{ and } \om^{\d}_t\in\R_+^d \mbox{ for all } t\in [0,1+\d] \big\}
	\e*
and the corresponding martingale optimal transport problem
\b*
\pp_{\d}&:=&\sup_{\P\in\Mc_+^{\d}}\E^{\P}[\z^{\d}(X^{\d})],
\e*
where similarly, $X^{\d}=(X^{\d}_t)_{0\le t\le k_{\d}}$ denotes the canonical process and $\Mc_+^{\d}$ denotes the set of martingale measures supported on $\Omo^{\d}$. The dual problem is slightly different. Denote further
 \b*
		\Omo^{c,\d}
		&:=&
		\big\{\om^{\d}\in\Omo^{\d}: \om^{\d}_{t^{\d}_i-}=\om^{\d}_{t^{\d}_i} \mbox{ for all } i=1,\cdots, m\big\}
	\e*
	and define the dual problem by
	\b*
		\dd^c_{\delta} 
		&:=& \inf_{(z, H) \in \Dc^c_{\delta}} z,
	\e*
with $\Dc^c_{\d}$ given by
\b*
		\Dc^c_{\d}
		&:=&
		\Big\{(z^{\d},H^{\d})\in \R\x\Ac^{\d}: z^{\d}+(H^{\d}\cdot \om^{\d})_1 \ge \zeta^{\d}(\om^{\d}) \mbox{ for all } \om^{\d}\in\Om_+^{c,\delta} \Big\},
		\e*
		where, similarly to Definition \ref{def:super-replication}, $\Ac^{\d}$ denotes the collection of all left-continuous adapted processes with bounded variation such that the stochastic integral $(H^{\d}\cdot X^{\d})$ is a supermartingale under all probability measures in $\Mc_+^{\d}$. 
		
	The main technical step for our result is the following.
		
	\begin{Lemma}\label{lemm:extension}
		Suppose that $\zeta$ satisfies \asss \ref{ass:xi_S-usc}, \ref{ass:xi_shift} and Condition \eqref{ass:xi_sko}. Then
		\be\label{ineq:invd}
			\dd^c_{{\d}}~\le~ \pp_{{\d}} \mbox{ for all } \d\ge 0.
		\ee
	\end{Lemma}
	
	The proof of Lemma \ref{lemm:extension} is adapted from Dolinsky \& Soner \cite{DS2} and is reported in Section \ref{sec:proof-lemma}.

	\begin{Lemma}\label{lemm:extension2}
		Suppose that $\zeta$ satisfies \ass \ref{ass:xi_S-usc}, \ref{ass:xi_shift} and Condition \eqref{ass:xi_sko}. Then
		\be\label{eq:asymp}
			\liminf_{\d \downarrow 0}\dd_{\d}^c~\ge~ \dd (\zeta)
			&\mbox{and}&
			\limsup_{\d \downarrow 0}\pp_{{\d}} ~\le~ \pp(\zeta).
		\ee
	\end{Lemma}
	
	\noindent \proof
	\rmi For each $(z^{\d}, H^{\d})\in \Dc^c_{{\d}}$ with $\d>0$ let us construct a robust super-replication on $\Omo$. For any $\om\in\Omo$ define $H_0(\om)=H^{\d}_0(\om^{\d,0})$ and 
	\b*
		H_t(\om)
		~:=~
		H^{\delta}_{t- t_i +t^{\d}_i} (\om^{\d,0}) \mbox{ for all } t\in (t_i, t_{i+1}] \mbox{ and } i=0,\cdots, m-1,
	\e*
where $ \om^{\d,0}\in\Omo^{c,\d}$ is defined as before by
\b*
\om^{\d,0}_t~=~\om_{(t-t_i^{\d}+t_i)\wedge t_{i+1}} \mbox{ for all } t\in [t_i^{\d}, t_{i+1}^{\d}] \mbox{ and } i=0,\cdots, m-1.
\e*
It is clear that $H$ is $\F-$adapted, left-continuous, with bounded variation, and $(H\cdot X)$ is a supermartingale under every $\P\in\Mco$, hence $H\in\Ac$. Moreover,
	\be \label{eq:zd_Hd}
		z^{\d}
		+
			\big(H^{\d}\cdot \om^{\d,0}\big)_{1+\d}
		~\ge~
			\zeta^{\delta} (\om^{\d,0}) \mbox{ for all } \om\in\Omo.
	\ee
	Notice that $\big(H^{\d}\cdot \om^{\d,0}\big)_{1+\d}=
		 (H \cdot\om)_1$, thus we obtain by \ass \ref{ass:xi_shift} and Condition \eqref{eq:zd_Hd}
		 \b*
		z^{\d}
		+
		(H\cdot \om)_1
		~\ge~
		\zeta(\om)-\et\Big(\frac{\d}{1+\d}\Big)\Big(1+\sum_{i=0}^m|\om_{t_i}|+\int_0^1|\om_t|dt\Big) \mbox{ for all } \om\in\Omo,
	\e*
	which yields $\dd_{\d}^c+\big(1+(m+2)d\big)\et\big(\frac{\d}{1+\d}\big) ~\ge~ \dd (\zeta)$ and therefore
	\b*
			\liminf_{\d \downarrow 0}\dd_{\d}^c&\ge& \dd (\zeta).
		\e*
 \rmii Let $(\d_n)_{n\ge 1}$ be such that $\delta_n > 0$ and $\delta_n \downarrow 0$.
	Then there is a sequence $(\P_n)_{n\ge 1}$ such that
	\b*
		\limsup_{n\to\infty}\pp_{\d_n}
		&=& \limsup_{n\to\infty}\E^{\P_n}\big[\zeta^{\delta_n}(X^{\d_n})\big].
	\e*
	For any fixed $\delta_0 > 0$, we assume without loss of generality that $\d_n\le\d_0$ for all $n\ge 1$.
	Then for each $n \ge 1$, let us define $\tilde{\P}_n:=\P_n\circ \big( \tilde{X}^{\d_n} \big)^{-1}$ 
	where $\tilde X^{\delta_n} (\om^{\d_n}) := X^{\d_0}(\om^{\d_0,\d_n})$ is the extended process from $\Om^{\delta_n}$ to $\Om^{\delta_0}$. It follows by Proposition \ref{prop:extension} \rmii that
	\b*
		\E^{\P_n} \big[ \zeta^{\delta_n}(X^{\d_n}) \big]
		&\le&
		\big(1+(m+2)d\big)\et\Big(\frac{\delta_0 - \delta_n}{1+\d_0}\Big) + \E^{\P_n} \big[ \zeta^{\delta_0}(\tilde X^{\delta_n}) \big] \\
	&=&
		\big(1+(m+2)d\big)\et\Big(\frac{\delta_0 - \delta_n}{1+\d_0}\Big) + \E^{\bar \P_n} \big[ \zeta^{\delta_0}(X^{\delta_0}) \big].
	\e*
Again by the same argument in Proposition \ref{prop:usc} we obtain
\b*
	\limsup_{n\to\infty}\pp_{\d_n}
		&\le&2\big(1+(m+2)d\big)\et\Big(\frac{\delta_0}{1+\d_0}\Big) +\pp(\z)
\e*
	which yields the required result since $\d_0>0$ is arbitrary. \qed

	\vspace{2mm}

	\noindent \textit{Proof of Theorem \ref{th:dwc}.} 
	Let $(z,H) \in \Dc(\zeta)$, we know by definition $z+(H\cdot \om)_1 \ge \zeta(\om)$, $\forall \om\in\Omo$.
	Taking expectation over each sides, it follows that
	\b*
		z~\ge~ \E^{\P}[\zeta(X)] \mbox{ for all } \P\in\Mco.
	\e*
	Then we get the weak duality $\pp(\zeta)\le \dd(\zeta)$.
	The reverse inequality follows by Lemmas \ref{lemm:extension} and \ref{lemm:extension2}. \qed


\subsection{Proof of Lemma \ref{lemm:extension}}\label{sec:proof-lemma}

The arguments are mainly adapted from Dolinsky \& Soner \cite{DS2} and the main idea is to discretize the paths on the Skorokhod space. By Proposition \ref{prop:extension} (i), the proof of $\dd_{\d}^c\le\pp_{\d}$ is not altered by the value of $\d$. We therefore consider $\d=0$ in this subsection.

\subsubsection{A probabilistic hedging problem}

For all $n \in\N $, put
	\b*
		A^{(n)}
		~:=~
		\big\{2^{-n}q: q \in\N^d\big\}
		&\mbox{and}&
		B^{(n)}
		~:=~
		\big\{i\sqrt{d}2^{-n}: i\in\N\big\}\cup\big\{\sqrt{d}2^{-n}/j: j\in\N^{\ast}\big\}.
	\e*
	We then define a subspace $\Omh:=\Omh^{(n)}\sbe\Omo$ as follows.

	\begin{Definition}\label{def:discret}
		A path $\om\in\Om_+$ belongs to $\Omh$ if there exist non-negative integers 
		$0=K_0<K_1+1<\cdots<K_{m}+m$ and a partition $\big\{0=\htt_0<\htt_1<\cdots<\htt_{K_m+m} =1\big\}$
		such that $\htt_{K_i+i}=t_i$ for $1\le i\le m$ and
		\b*
			\om_t 
			&=&
			\sum_{i=0}^{m-1}
			\Big(
				\sum_{k=K_i+i}^{K_{i+1}+i-1} \om_{{\htt_k}}  \mathds 1_{[\htt_k,\htt_{k+1})}(t)
				+
				\om_{t_{i+1}}  \mathds 1_{[\htt_{K_{i+1}+i},t_{i+1})}(t)\Big)
			+
			\om_1 \mathds 1_{\{t=1\}},
		\e*
		where $\om_{t_i}\in A^{(n)}$ for $1\le i\le m$ and for $0\le i< m$
		\b*
			&&\om_{\htt_k}~\in~ A^{(n+k-K_i-i)},~ K_i+i< k<K_{i+1}+i+1, \\
			&&\htt_k-\htt_{k-1}~\in~ B^{(n+k-K_i-i)},~ K_i+i<k<K_{i+1}+i+1.
		\e*
	\end{Definition}

	 Notice that $\Omh$ is countable, then there exists a probability measure $\hP:=\hP^{(n)}$ on $\Omo$ supported on $\Omh$ which gives positive weight to every element of $\Omh$. 
	In particular, the canonical process $X$ has finitely many jumps $\hP$-almost surely.
	Denote by $\hF$ the completed filtration of $\F$ under $\hP$. 
	Put
	\b*
		\hH^{(n)}
		&:=&
		\Big\{ 
			\hh: [0,1]\x\Omo\to\R^d \mbox{ is } \hF-\mbox{predictable such that } \|\hh\|\le n 
		\Big\}
	\e*
	and
	\b*
		\hA^{(n)}
		&:=&
		\Big\{
			\hh \in \hH^{(n)}:
			(\hh\cdot X)_t \ge K \mbox{ for all } t\in [0,1],~ \hP\mbox{-a.s. for some } K\in\R
		\Big\}.
	\e*
	Let
	\b*
		\hD^{(n)} (\zeta)
		&:=&
		\Big\{(z,\hh)\in \R\x\hA^{(n)}:z+(\hh\cdot X)_{1}~\ge~\zeta(X),~ \hP\mbox{-a.s.}\Big\}
	\e*
	and define the robust superhedging problem under the dominating measure $\hP$
	\b*
		\dd^{(n)}(\zeta)&:=&\inf_{(z,\hh) \in \hD^{(n)}(\zeta)} z.
	\e*
Let $\hat{\Pc}\sbe\Pc$ be the subset of probability measures supported on $\Omh$,
	and $\hat{\Mc}_n\sbe\hat{\Pc}$ be the subset of probability measures $\Q$ that have the following properties:
	\b*
		\E^{\Q}\Big[\sum_{k=1}^{K_m+m}\Big|\E^{\Q}[X_{\htt_k}|\Fc_{\htt_k-}]-X_{\htt_{k-1}}\Big|\Big]&\le&\frac{1}{n},
	\e*
	where $0<\htt_1(\om)<\cdots<\htt_{K_m+m-1}(\om)<1$ are the jumps times of the piecewise constant process $X(\om)$ with $\htt_0(\om)=0$ and $\htt_{K_m+m}(\om)=1$. Then the required result $\dd^c(\zeta)\le\pp(\zeta)$ follows from the following Propositions \ref{prop:approxp} and \ref{prop:approxd}.

	\begin{Proposition}\label{prop:approxp}
		Assume that $\z$ satisfies Assumptions \ref{ass:xi_S-usc}, \ref{ass:xi_shift} and is $\L^{\infty}-$uniformly continuous, then
		\b*
			\limsup_{n\to\infty} \dd^{(n)}(\zeta)
			&\le&
			\pp(\zeta).
		\e*
	\end{Proposition}
	
\noindent \proof	 \rmi 	From Example 2.3 and Proposition 4.1 in F\"ollmer \& Kramkov \cite{FK},
	it follows that
	\b*
		\dd^{(n)}(\zeta)
		&=&
		\sup_{\Q\in\hat{\Pc}} 
		\E^{\Q} \Big[
			\zeta-n\sum_{k=1}^{K_m+m}\Big|X_{\htt_{k-1}}-\E^{\Q}[X_{\htt_k}|\Fc_{\htt_{k}-}]\Big|
		\Big].
	\e*
	Since $0\le \zeta\le 1$, we determine that $\dd^{(n)}(\zeta)\ge 0$ 
	and we have for every $\Q\in\hat{\Pc}\backslash\hat{\Mc}_n$,
	\b*
		\E^{\Q}\Big[\zeta-n\sum_{k=1}^{K_m+m}\Big|X_{\htt_{k-1}}-\E^{\Q}[X_{\htt_k}|\Fc_{\htt_{k}-}]\Big|\Big]
		&\le& 0,
	\e*
	which yields 
	\b*
			\dd^{(n)}(\zeta)&\le&\sup_{\Q\in\hat{\Mc}_n}\E^{\Q}[\zeta(X)].
		\e*
\rmii Let us take a sequence $(\Q_n)_{n\ge 1}$ with $\Q_n\in\hat{\Mc}_n$ such that
	\b*
		\limsup_{n\to\infty}\big\{\sup_{\Q\in\hat{\Mc}_n}\E^{\Q}[\zeta(X)]\big\}
		&=&
		\limsup_{n\to\infty}\E^{\Q_n}[\zeta(X)].
	\e*
	Since under each $\Q_n$ the canonical process $X$ is piecewise constant 
	with jump times $0<\htt_1<\cdots<\htt_{K_m+m-1}<1$, 
	$X$ is a $\Q_n-$semimartingale.
	Then we have the decomposition $X=M^{\Q_n}-A^{\Q_n}$, 
	where $A^{\Q_n}$ is a predictable process of bounded variation and $M^{\Q_n}$ is a martingale under $\Q_n$. Moreover, $A^{\Q_n}$ is identified by
	\b*
		A_t^{\Q_n}~=~\sum_{k=1}^{K_m+m-1}\mathds{1}_{[\htt_k,\htt_{k+1})}(t)\sum_{j=1}^k\Big[X_{\htt_{j-1}}-\E^{\Q_n}[X_{\htt_j}|\Fc_{\htt_j-}]\Big] \mbox{ for all } t\in [0,1),
	\e*
	and $A_{1}^{\Q_n}=\lim_{t\to 1}A_t^{\Q_n}$.
	It follows then
	$\E^{\Q_n}\big[|X_{1}-M_{1}^{\Q_n}|\big]\le\E^{\Q_n}\big[|A^{\Q_n}_{1}|\big]\le 1/n$
	and
	\b*
		\Q_n\Big[\|A^{\Q_n}\|\ge n^{-1/2}\Big]
		~~\le~~
		n^{1/2}\E^{\Q_n}\Big[\sum_{k=1}^{K_m+m-1}\big|X_{\htt_{k-1}}-\E^{\Q_n}[X_{\htt_k}|\Fc_{\htt_k-}]\big|\Big]
		~~\le~~ 
		n^{-1/2}.
	\e*
	Since $\zeta$ is $\L^{\infty}-$uniformly continuous, one obtains
\b*
		\limsup_{n\to\infty}\E^{\Q_n}[\zeta(X)]
		&\le& 
		\limsup_{n\to\infty}\E^{\Q_n}[\zeta(M^{\Q_n})].
	\e*
	Let $\P_n=\Q_n\circ (M^{\Q_n})^{-1}$, then
	\b*
		\sup_{n\ge 1}\E^{\P_n}[|X_{1}|] 
		&=&
		\sup_{n\ge 1}\E^{\Q_n}[|M_{1}^{\Q_n}|] \\
			&\le&
			\sup_{n\ge 1}\E^{\Q_n}[|M_{1}^{\Q_n}-X_{1}|]
			+
			\sup_{n\ge 1}\E^{\Q_n}[X_{1}]\\
		&\le&
		\sup_{n\ge 1}\E^{\Q_n}[|M_{1}^{\Q_n}-X_{1}|]+ \sup_{n\ge 1}\E^{\Q_n}[X_{1}-M_{1}^{\Q_n}]+\sup_{n\ge 1}\E^{\Q_n}[M_{1}^{\Q_n}]\\
		&\le &
		1+\frac{2}{n}
		~~\le~~ 3.
	\e*
	By Assumptions \ref{ass:xi_S-usc} and \ref{ass:xi_shift}, it follows that for any $\bep\in\R^m_+$ such that $0<|\bep|<\Delta \T$
	\b*
		\limsup_{n\to\infty}\E^{\P_n}[\z(X)]
		&\le&
		\limsup_{n\to\infty}\E^{\P_n}[\z(X_{f_{\bep}})]+\big(1+(m+2)d\big)\a(|\bep|).
	\e*
	Again with the same reasoning, we may prove 
	\b*
		\limsup_{n\to\infty}\E^{\P_n}[\z(X_{f_{\bep}})]&\le&\pp(\z).
	\e*
	Since $\bep$ is arbitrary we get
	\b*
		\limsup_{n\to\infty}\Big[\sup_{\Q\in\hat{\Mc}_n}\E^{\Q}[\zeta(X)]\Big]&\le& \pp(\zeta),
	\e*
	and hence the required result. \qed

\subsubsection{Time-space discretization}
	
	\noindent\underline{\it Discretization}: 
	For each $\om\in\Omo^c$ let us define $\t_k:=\t^{(n)}_k(\om)$ and $K_i:= K^{(n)}_i(\om)$ by
	\b*
		\t_0&:=&0,
		~~~~K_0~:=~0, \\
		\t_1 
		&:=&
		t_1\wedge\sqrt{d}2^{-n}\wedge\inf\crl{t>0: |\om_t-\om_0|\ge 2^{-n}}, \\
		\tau_{k+1}
		&:=&
		t_1\wedge\brak{\t_k+\Dl\t_k}\wedge\inf\crl{t>\t_k: |\om_t-\om_{\t_k}|\ge 2^{-n}},
		~ \Dl\t_k=\t_k-\t_{k-1} \mbox{ for } k\ge 1.
	\e*
	Set further
	\b*
		K_1
		&:=&
		\min\crl{k\in\N: \t_k=t_1}.
	\e*
	Recursively, we define for $1\le i\le m-1$ and $k\ge K_i$,
	\b*
		\t_{K_i+1}&:=&t_{i+1}\wedge(t_i+\sqrt{d}2^{-n})\wedge\inf\crl{t>t_i: |\om_t-\om_{t_i}|\ge 2^{-n}}, \\
		\t_{k+1}&:=&t_{i+1}\wedge\brak{\t_k+\Dl\t_k}\wedge\inf\crl{t>\t_k: |\om_{t}-\om_{\t_k}|\ge 2^{-n}} \mbox{ for } k\ge K_i+1
	\e*
	and
	\b*
		K_{i+1}&:=&\min\crl{k\in\N: \t_k=t_{i+1}}.
	\e*
	 Notice that the above $\tau_k$, $k \ge 0$ are all stopping times w.r.t. to the right-continuous filtration $\F^+ = (\Fc_{t+})_{t \ge 0}$,
	and
	\b*
		0~=~\t_0~<~\t_1~\cdots~<~\t_{K_m}~=~1
		&\mbox{and}&
		\t_{K_i}~=~t_i \mbox{ for all } i=1,\cdots, m.
	\e*
	Moreover, for $0\le i\le m-1$, $K_i< k\le K_{i+1}$ and $ t\in [\t_{k-1},\t_{k})$,
	\b*
		|\om_t-\om_{\t_{k-1}}|~\le~ 2^{-n} &\mbox{and}&
		\Dl\t_{k+1}~\le~ \Dl\t_k~\le~ 2^{-n}.
	\e*
	Also by the continuity of $\om$ at $\t_{K_i}=t_i$ for all $i=1, \cdots, m$
	\b*
		|\om_t-\om_{\t_{K_i-1}}|~\le~ 2^{-n} \mbox{ for all } t\in [\t_{K_{i}-1},t_i] \mbox{ and } i=1,\cdots, m.
	\e*

\noindent\underline{\it Lifting}:
	Set $\htt_0:=0$ and for $0\le i\le m-1$
	\b*
		\htt_{K_i+1}&:=&\htt_{K_i}+\sqrt{d}2^{-n}, \\
		\htt_k&:=&\htt_{k-1}+\big(1-\sqrt{d}2^{-n}/\Dl t_{i+1}\big)\sup\big\{\Dl t>0: \Dl t\in B^{(n+k-K_i-i)},~ \Dl t<\Dl\t_{k-1}\big\}, \\
		&&\mbox{for all } K_i+i+2\le k\le K_{i+1}+i, \\
		\htt_{K_{i+1}+i+1}&:=&t_{i+1}.
	\e*
	Denote  $\hp(\om)=\big(\hp(\om)\big)_{0\le t\le 1}$ by
		\b*
		\hp_t(\om)
		&:=&
		\sum_{i=0}^{m-1}\Big\{\sum_{k=K_i+i}^{K_{i+1}+i-1}\pi^{(n+k-K_i-i)}(\om_{\t_k}) \mathds{1}_{[\htt_k,\htt_{k+1})}(t)+\pi^{(n)}(\om_{t_i})\mathds{1}_{[\htt_{K_r},t_{i+1})}(t)\Big\} \\
		&& +
		\pi^{(n)}(\om_1)\mathds{1}_{\{t=1\}},
	\e*
	then $\hp(\om)\in\Omh$. 
	For each $\hh\in\hA$ we may define
	\be \label{eq:star}
		H_t(\om)~:=~\sum_{k=0}^{K_{m}-1}\hh_{\htt_{k+1}(\om)}(\hp(\om))\mathds{1}_{(\t_k(\om),\t_{k+1}(\om)]}(t) \mbox{ for all } (\om, t)\in \Omo\x [0,1].
	\ee
	We observe that $(\om,t) \mapsto H_t(\om)$ is Borel measurable on $\Om_+ \x [0,1]$,
	$t \mapsto H_t(\om)$ is left-continuous,
	$\om \mapsto H_t(\om)$ is $\Fc_{t+}-$measurable.
	Hence $H$ is $\F^+-$predictable, which is equivalent to be $\F-$predictable.
	Further, following the argument of \lems 3.5 and 3.6 of Dolinsky \& Soner \cite{DS2}, we see that the process $H$ defined by \eqref{eq:star} belongs to $\Ac$, and more importantly, there exists some constant  $C>0$ independent of $n$ such that for all $\om\in\Omo$,
		\be\label{ineq:shift}
			\r_{\T}(\om,\hp(\om))\le C2^{-n}\big(1+\|\om\|\big)
			\mbox{ and }
			\big|(H \cdot \om)_{1}-(\hh(\hp(\om))\cdot \hp(\om))_{1}\big|
			\le
			Cn2^{-n}.
		\ee

	\begin{Proposition}\label{prop:approxd}
		Assume that $\z$ satisfies Condition \eqref{ass:xi_sko}, then one has
		\b*
			\liminf_{n\to\infty}\dd^{(n)}(\zeta)&\ge& \dd^c(\zeta).
		\e*
	\end{Proposition}
	
	\noindent \proof Take an arbitrary $(\hat z,\hh)\in\hD$. 
	Then for any $\om\in\Omo$ one has $\hp(\om)\in\Omh$ and thus
	\b*
		\hat z+(\hh(\hp(\om))\cdot\hp(\om))_{1}~\ge~\z(\hp(\om)) \mbox{ for all } \om\in\Omo.
	\e*
	Take $H$ constructed as \eqref{eq:star}, then by \eqref{ineq:shift}, we have $H\in\Ac$ and 
	\b*
		\hat z+(H\cdot \om)_{1}~\ge~\zeta(\hp(\om))-Cn2^{-n}  \mbox{ for all } \om\in\Omo.
	\e*
	Moreover, by the construction of $\hp(\om)$ one has $\|\hp(\om)\|\le\|\om\|$. Notice that $\z$ is boundedly supported, saying by $\big\{\om\in\Om: \|\om\|\le R\big\}$. Then by \eqref{ass:xi_sko} one has a continuous increasing function $\bt: \R_+\to\R_+$ with $\bt(0)=0$ such that for all $\om\in\Omo$,
	\b*
		\zeta(\hp(\om))
		~~\ge~~
		\zeta(\om)-\mathds{1}_{\{\|\om\|\le R\}}\bt \big(\r_{\T}\big(\om,\hp(\om)\big)\big)
		~~\ge~~
		\zeta(\om)-\bt\big(C(1+R)2^{-n}\big),
	\e*
	which implies that $\big(\hat z+\bt\big(C(1+R)2^{-n}\big)+Cn2^{-n}, H\big)\in\Dc^c(\z)$. Hence
	\b*
		\dd^c(\zeta)&\le& \dd^{(n)}(\zeta)+\bt\big(C(1+R)2^{-n}\big)+Cn2^{-n},
	\e*
	which yields the required result. \qed

\subsection{Proof of Proposition \ref{prop:extension} (i)} \label{subsec:extension}

Recal that $\xi$ satisfies \asss \ref{ass:xi_S-usc}, \ref{ass:xi_shift} and Condition \eqref{ass:xi_sko}. The required statement follows from the three lemmas below.

\begin{Lemma}\label{lem:ex1}
	$\limsup_{n\to\infty}\xi^{\d}(\om^{\d,n})\le\xi(\om^{\d,0})$ holds for any sequence $(\om^{\d, n})_{n\ge 0}\sbe\Om^{\d}$ such that 
	\b*
	\om^{\d, n}~\sto~\om^{\d, 0} &\mbox{and}& \om^{\d, n}_{t_i^{\d}}~\lro\om^{\d, 0}_{t_i^{\d}} \mbox{ for all } i=0,\cdots ,m-1.
	\e*
\end{Lemma}
	
\noindent \proof For the sake of simplicity we may assume that
	\b*
	\limsup_{n\to\infty}\xi^{\d}(\om^{\d, n})&=&\lim_{n\to\infty}\xi^{\d}(\om^{\d, n}).
	\e* 
	Since $(\om^{\d, n})_{n\ge 1}$ is $S-$tight, then by the $S-$tightness criteria and the construction in \eqref{def:extension} we determine that $(\bar \om^{\d, n})_{n\ge 1}$ is again $S-$tight, which yields a convergent subsequence $(\bar \om^{\d, n_{k}})_{k\ge 1}$ and a limit $\om^0\in\Om$, i.e. $\bar \om^{\d, n_{k}}\sto\om^0$. Clearly,
$\om^{\d,n}_{t_i^{\d}}\lro\om^{\d,0}_{t_i^{\d}}$ implies in particular that $\bar\om^{\d,n_{k}}_{t_i}\lro\bar\om^{\d,0}_{t_i}$ for all $i=0,\cdots, m-1$. Next, $\om^{\d, n_{k}}\sto\om^{\d,0}$ yields a countable set $\Tc\sbe[0,1+\d)$ such that 
\b*
\om^{\d, n_{k}}_t~\lro~\om^{\d,0}_t \mbox{ for all } t\in [0,1+\d]\backslash\Tc,
\e*
which yields another countable set $\Tc'\sbe [0,1)$ such that
\b*
\bar \om^{\d, n_{k}}_t~\lro~\bar \om^{\d,0}_t \mbox{ for all } t\in [0,1]\backslash\Tc'.
\e*	
Hence one has $\bar \om^{\d,0}=\om^0$ and thus
\b*
	\bar \om^{\d, n_k}~\sto~\bar \om^{\d, 0} &\mbox{and}& \bar\om^{\d, n_k}_{t_i}~\lro~\bar \om^{\d, 0}_{t_i} \mbox{ for all } i=0,\cdots, m-1,
	\e*
which implies that
\b*
\lim_{k\to\infty}\xi^{\d}(\om^{\d, n_{k}})~~=~~\lim_{k\to\infty}\xi(\bar \om^{\d, n_{k}})&\le&\xi(\bar \om^{\d, 0})~~=~~\xi^{\d}( \om^{\d, 0}).
\e*	
\qed
	
	\begin{Lemma}\label{lem:ex2}
		There exists a continuous function  $\a_{\d} : \R_+ \to \R_+$ with $\a_{\d}(0)=0$ such that for all $\bep=(\eps_1, \cdots, \eps_m)\in \R_+^m$ sufficiently small one has
		\b*
			\big|\xi^{\d}(\om^{\d}) -\xi^{\d}\big(\om^{\d}_{ f^{\delta}_{\bep}} \big)\big|,~~ \big|\xi^{\d}(\om^{\d}) -\xi^{\d}\big(\om^{\d}_{ b^{\delta}_{\bep}} \big)\big|
			&\le&
			\a_{\d}(|\bep|)\Big(1+\sum_{i=0}^m|\om^{\d}_{t_i^{\d}}|+\int_0^{1+\d}|\om^{\d}_t|dt\Big), 
		\e* 
		where $f_{\bep}^{\d},~ b_{\bep}^{\d} : [0, 1+ \delta] \to [0, 1+ \delta]$ are two non-decreasing functions defined as in \eqref{forward} and \eqref{backward}.
	\end{Lemma}
	
	\noindent \proof 
	We only prove the inequality on $f^{\d}_{\bep}$, while the inequality on $b^{\delta}_{\bep}$ follows by the same arguments.
Define $\bar f_{\bep}^{\d}: [0,1+\d]\to [0,1+\d]$ by
\b*
\bar f_{\bep}^{\d}(t)~:=~\frac{1}{k_{\d}}f_{\bep}^{\d}(k_{\d}t) \mbox{ for all } t\in [0,1].
\e*
Thus by the construction of $f_{\bep}^{\d}$ we get
\b*
\bar f_{\bep}^{\d}(t)&=&\frac{1}{k_{\d}}\sum_{i=1}^m
			\mathds{1}_{(t_{i-1}^{\d}, t_{i}^{\d}]}(k_{\d}t)\Big(t_{i-1}^{\d}+\frac{\Dl t_i^{\d}}{\Dl t_i^{\d}-\eps_i}\big(k_{\d}t-t^{\d}_{i-1}-\eps_i\big)^+\Big) \\
			&=&\sum_{i=1}^m
			\mathds{1}_{(t_{i-1}, t_{i}]}(t)\Big(t_{i-1}+\frac{\Dl t_i}{\Dl t_i-\eps_i /k_{\d}}\big(t-t_{i-1}-\frac{\eps_i}{k_{\d}}\big)^+\Big), 
\e*
	which implies that
	\b*
	\bar \om^{\d}_{f^{\d}_{\bep}(t)}~=~\om^{\d}_{k_{\d}\bar f_{\bep}^{\d}(t)} \mbox{ for all } t\in [0,1]
	\e*
	and thus $\bar \om^{\d}_{f^{\d}_{\bep}}=\bar \om^{\d}\circ \bar f_{\bep}^{\d}$. Hence
	\b*
\big|\xi^{\d}(\om^{\d}) -\xi^{\d}\big(\bar \om^{\d}_{f^{\d}_{\bep}}\big)\big|
&=&\big|\xi(\bar \om^{\d}) -\xi\big(\bar \om^{\d}\circ \bar f_{\bep}^{\d} \big)\big| \\
&\le &\a(|\bep|/ k_{\d})\Big(1+\sum_{i=0}^m|\bar \om^{\d}_{t_i}|+\int_0^{1}|\bar \om^{\d}_t|dt\Big) \\
&= &\a(|\bep|/ k_{\d})\Big(1+\sum_{i=0}^m|\om^{\d}_{t_i^{\d}}|+\frac{1}{k_{\d}}\int_0^{1+\d}|\om^{\d}_t|dt\Big) \\
&\le &\a(|\bep|/ k_{\d})\Big(1+\sum_{i=0}^m|\om^{\d}_{t_i^{\d}}|+\int_0^{1+\d}|\om^{\d}_t|dt\Big). 
\e*
The proof is completed by taking $\a_{\d}(\cdot)=\a(\cdot /k_{\d})$. \qed

\begin{Lemma}\label{lem:ex3}
$\xi^{\d}$ is $L^{\infty}-$uniformly continuous and satisfies Condition \eqref{ass:xi_sko} for $\r_{\T^{\d}}$.
\end{Lemma}

\noindent \proof For any $\om^{\d}, v^{\d}\in\Om^{\d}$ such that $\|v^{\d}\|\le\|\om^{\d}\|$, one has 
	\b*
	\xi^{\d}(\om^{\d})-\xi^{\d}(v^{\d})&=&\xi(\bar \om^{\d})-\xi(\bar v^{\d})~~\le~~\bt\big(\r_{\Tf}(\bar \om^{\d},\bar v^{\d})\big), \\
	\big|\xi^{\d}(\om^{\d})-\xi^{\d}(v^{\d})\big|&=&\big|\xi(\bar \om^{\d})-\xi(\bar v^{\d})\big|.
	\e*
It is thus enough to show that
	\b*
	\r_{[t_{i-1},t_i]}(\bar \om^{\d},\bar v^{\d})&\le& \r_{[t_{i-1}^{\d},t_i^{\d}]}(\om^{\d},v^{\d}) \mbox{ for all } i=1, \cdots, m
	\e*
and
\b*
\Big|\int_0^{1}\big(\bar\om^{\d}_t-\bar v^{\d}_t\big)dt\Big|&\le&\Big|\int_0^{1+{\d}}\big(\om^{\d}_t-v^{\d}_t\big)dt\Big|.
\e*
Let $\Gm_{[s,t]}$ denotes the collection of strictly increasing continuous functions $\g$ defined on $[s,t]$ such that $\g(s)=s$ and $\g(t)=t$. For any $\g^{\d}\in\Gm_{[t_{i-1}^{\d},t_i^{\d}]}$, define $\g\in\Gm_{[t_{i-1},t_i]}$ by
\b*
\g(t)~:=~\frac{1}{k_{\d}}\g^{\d}\big(k_{\d}t\big) \mbox{ for all } t\in [t_{i-1},t_i].
\e* 
Hence
\b*
&&\sup_{t_{i-1}\le t\le t_i}\big|\bar\om^{\d}_ {\g(t)}-\bar v^{\d}_t\big|
~~=~~\sup_{t_{i-1}\le t\le t_i}\big|\om^{\d}_{k_{\d}\g(t)}-v^{\d}_{k_{\d}t}\big|
~~=~~\sup_{t_{i-1}\le t\le t_i}\big|\om^{\d}_{\g^{\d}(k_{\d}t)}-v^{\d}_{k_{\d}t}\big| \\
&=&\sup_{t_{i-1}^{\d}\le t\le t_i^{\d}}\big|\om^{\d}_{\g^{\d}(t)}-v^{\d}_t\big|
\e*
and
\b*
&&\sup_{t_{i-1}\le t\le t_i}\big|\g(t)-t\big|~~=~~\sup_{t_{i-1}\le t\le t_i}\Big|\frac{1}{k_{\d}}\g^{\d}\big(k_{\d}t\big)-t\Big|
~~=~~\frac{1}{k_{\d}}\sup_{t_{i-1}\le t\le t_i}\big|\g^{\d}\big(k_{\d}t\big)-k_{\d}t\big| \\
&=&\frac{1}{k_{\d}}\sup_{t_{i-1}^{\d}\le t\le t_i^{\d}}\big|\g^{\d}(t)-t\big|
~~\le~~\sup_{t_{i-1}^{\d}\le t\le t_i^{\d}}\big|\g^{\d}(t)-t\big|,
\e*
which implies that
\b*
\r_{[t_{i-1},t_i]}(\bar\om^{\d},\bar v^{\d})&\le&\r_{[t_{i-1}^{\d},t_i^{\d}]}(\om^{\d}, v^{\d}).
\e*
We may thus conclude by
\b*
\Big|\int_0^{1}\big(\bar\om^{\d}_t-\bar v^{\d}_t\big)dt\Big|
&=&\frac{1}{k_{\d}}\Big|\int_0^{1+\d}\big(\om^{\d}_t-v^{\d}_t\big)dt\Big|
~~\le~~ \Big|\int_0^{1+\d}\big(\om^{\d}_t-v^{\d}_t\big)dt\Big|.
\e*
\qed


\appendix

\section{Appendix}

\subsection{Tightness under $S-$topology}

	Recall that $\Om= \D([0,1],\R^d)$ is the Skorokhod space of c\`adl\`ag paths on $[0,1]$, with canonical process $X=(X_t)_{0\le t\le 1}$ and canonical filtration $\F = (\Fc_t)_{0 \le t \le 1}$, and $\Pc$ denotes the set of all probability measures on $(\Om, \Fc_1)$. A sequence of probability measures $(\P_n)_{n\ge 1}\subset\Pc$ is said to be $S-$tight if for any $\eps>0$, there exists a $S-$compact set $K_{\eps}\subset\Om$ such that 
	\b*
	\inf_{n\ge 1}\P_n\big[X\in K_{\eps}\big]&\ge&1-\eps.
	\e*
	The following result is recalled from Jakubowski \cite{Jakubowski} (see their Theorem 3.11 and the discussion at the beginning of Section 4 in \cite{Jakubowski}).
	
		\begin{Theorem}[Jakubowski]\label{Thm:Jakubowski}
		\rmi Let $(\P_n)_{n\ge 1}\sbe\Pc$ be a sequence of probability measures such that $X$ is a $\P_n-$supermartingale for all $n\ge 1$, then
		\b*
			\sup_{n\ge 1}\sup_{0\le t\le 1}\E^{\P_n} \big[ |X_t| \big]
			~<~
			+\infty
			~~&\Longrightarrow&~~
			(\P_n)_{n\ge 1} ~\mbox{is}~ S-\mbox{tight}.
		\e*
\rmii Let $(\P_n)_{n\ge 1}\sbe\Pc$ be a $S-$tight sequence of probability measures.
		Then there exist a subsequence $(\P_{n_k})_{k\ge 1}$, 
		a probability measure $\P\in\Pc$  and a countable subset $\Tc\sb [0,1)$ 
		such that for all finite sets $\{u_1<u_2<\cdots<u_r\}\sb [0,1]\backslash\Tc$,
		\be \label{eq:cvg_marginals}
			\P_{n_k}\circ (X_{u_1}, \cdots,  X_{u_r})^{-1}
			~\lro~
			\P\circ (X_{u_1}, \cdots, X_{u_r})^{-1}
			\mbox{ as } k\to\infty.
		\ee
		In particular, $X^{n_k}\dwo X^0$ as $k\to\infty$.
\end{Theorem}

\subsection{Dual space of $\Mf$}\label{lem:dualspace}

	Recall that $\Mf$ denotes the space of all finite signed measures $\m$ on $\R^d$ satisfying
	\b*
		\int_{\R^d}\big(1+|x|\big)|\m|(dx)&<&+\infty,
	\e*
	and it is equipped with the topology induced by the convergence $\woi$.
	We would like identify its dual space as $\Lambda$,
	where the arguments are mainly adapted from Lemma 3.2.3 of Deuschel \& Stroock \cite{DeuschelStroock}.
	Notice that the topology on $\Mf$ is generated by all the following open balls
	\b*
		U_{\lambda^1, \cdots, \lambda^m , c}(\mu)
		&:=&
		\Big\{ \nu \in\Mf: \big|\mu (\lambda^i)-\nu (\lambda^i)\big|
			<c \mbox{ for all } 1\le i \le m 
		\Big\},
	\e*
	where $\lambda^i \in \Lb$ for $1\le i\le m$ and $c>0$.
	Let $\Oc$ be the collection of open sets generated by the open balls above, 
	then clearly, every open set $U\in\Oc$ could be expressed as
	\b*
		U
		~=~
		\bigcup_{\alpha} U_{\lambda^1_{\a}, \cdots, \lambda_{\a}^{n_{\a}}, c^{\a}}(\mu^{\a}) 
		~\mbox{with}~
		\lambda^i_{\alpha} \in \Lb \mbox{ for } 1\le i\le n_{\a}, n_{\a}\in\N \mbox{ and } c^{\a}>0.
	\e*

	\begin{Theorem}
		The space $(\Mf, \Oc)$ is a Hausdorff locally convex space,
		whose dual space can be identified by $\Mf^{\ast}=\Lb$.
	\end{Theorem}

	\noindent \proof \rmi First, $(\Mf, \Oc)$ is clearly a topological vector space.
	For every $\mu \in\Mf$, let
	\b*
		\Uc(\mu)
		&:=&
		\big\{ U_{\lambda^1_{\a}, \cdots, \lambda_{\a}^{n_{\a}}, c^{\a}}(\mu)
			: 
			\lambda^i_{\a}\in\Lb_1 
			\mbox{ for }1\le i\le n_{\a}, n_{\a}\in\N \mbox{ and } c^{\a}>0
		\big\}.
	\e*
	By definition, one can check that $\Uc(\mu)$ is a local basis of $\mu$ for every $\mu \in \Mf$.
	Moreover, by denoting $\0 \in \Mf$ the null measure,
	$\Uc( \0 )$ is a local basis of absolutely convex absorbent sets and thus $\Mf$ is a locally convex space.

	\vspace{1mm}

	\noindent \rmii Now, let us identify the dual space of $\Mf$.
	First, for every $\lambda \in \Lb$, the map $F_{\lambda} : \Mf \to \R$ defined by 
	$F_{\lambda}( \mu) := \mu(\lambda)$ gives a unique element in $\Mf^{\ast}$,
	and hence $\Lb \sbe \Mf^{\ast}$. On the other hand,  for any $F\in \Mf^{\ast}$, we define a function $\l^F$ by
	\b*
		\l^F(x)~:=~F(\d_{\{x\}}) \mbox{ for all } x\in\R^d.
	\e*
	Clearly one has the following implication
	\b*
		x_n\to x_0
		~~\Longrightarrow~~
		\d_{\{x_n\}}\stackrel{\Oc}{\lro}\d_{\{x_0\}}
		~~\Longrightarrow~~
		\l^F(x_n)\to\l^F(x_0),
	\e*
	which implies that $\l^F$ is continuous.
	It follows that the set $F^{-1}\big((-1,1)\big)$ is open and thus there exists some 
	$U_{\l^1, \cdots, \l^m, c}(\mathbf 0)$ such that
	\b*
		U_{\l^1, \cdots, \l^m, c}(\mathbf 0)&\sbe&F^{-1}\big((-1,1)\big),
	\e*
	where $\l^i\in\Lb$ for all $i=1,\cdots, m$ and $c>0$. 
	Now for any $\m\in\Mf$ such that $\sum_{i=1}^{m}\big|\m(\l^i)| >  0$,
	we define
	\b*
		\bar \m&:=&\frac{c\m}{\sum_{i=1}^{m}\big|\m(\l^i)\big|}.
	\e*
	Then $\bar \m\in U_{\l^1, \cdots, \l^m, c}(\mathbf 0)$ and thus $|F(\bar \m)|<1$. It follows that
	\b*
		|F(\m)|
		~\le~
		c\sum_{i=1}^{m}\big|\m(\l^i)| \mbox{ for all } \m\in\Mf,
	\e*
	and hence $\l^F\in\Lb$. 
	When $\mu$ is a linear combination of Dirac measures, it is obvious that $F(\m)=\m(\l^F)$.
	Moreover, since such $\mu$ are dense in $\Mf$, it follows that $F(\m)=\m(\l^F)$ holds for all $\mu \in \Mf$. \qed


\end{document}